\documentclass[reqno]{amsart}

\usepackage{amsmath}
\usepackage{amssymb}
\usepackage{amsthm}
\usepackage{xspace}
\usepackage{pictexwd}
\usepackage{float}  \restylefloat{figure} %\restylefloat{table}
\usepackage{amscd}
\usepackage{enumerate}

\theoremstyle{plain}
\newtheorem{thm}{Theorem}[section]
\newtheorem*{thm*}{Theorem}
\newtheorem{prop}[thm]{Proposition}
\newtheorem{lem}[thm]{Lemma}
\newtheorem{cor}[thm]{Corollary}

\theoremstyle{definition}
\newtheorem{mydef}{Definition}[section]

\theoremstyle{remark}
\newtheorem*{rem}{Remark}
\newtheorem*{notation}{Notation}

\def\forcehmode{\hskip0pt\relax}

\let\myskip=\medskip

\def\definebb#1=#2.{\def#1{{{\mathbb #2}^{\vphantom{x}}}}}
\definebb \c=C.
\definebb \r=R.
\definebb \s=S.
\definebb \z=Z.
\definebb \h=H.
\def\Cal#1{\mathcal{#1}}

\def\calA{{\Cal A}}

\def\num{n}

\def\tMod{\widetilde{\Mod}}
\def\tTgll{\tilde T_{g,l_h,l_p}}
\def\Tgll{T_{g,l_h,l_p}}
\def\Mgll{M_{g,l_h,l_p}}
\def\Modgll{\Mod_{g,l_h,l_p}}
\def\Smgll{S^m_{g,l_h,l_p}}
\def\Modmgll{\Mod^m_{g,l_h,l_p}}

\def\lev{{s}}
\def\levm{\lev_m}
\def\lat{\Ga}

\def\lats{\lat^*}

\def\tG{{\tilde G}}
\def\Gm{{G_m}}

\def\bA{{\bar A}}
\def\bB{{\bar B}}
\def\bC{{\bar C}}
\def\bD{{\bar D}}
\def\bV{{\bar V}}

\def\bX{{\bar X}}

\def\zd{\z_2}
\def\zm{\z_m}

\def\abcd{{a\,b\choose c\,d}}

\def\ABCD{\begin{pmatrix} a&b \\ c&d \end{pmatrix}}

\def\hyp{\h}

\def\PSL{\MathOpPSL(2,\r)}
\def\SL{\MathOpSL(2,\r)}

\def\tPSL{\widetilde{\MathOpPSL}(2,\r)}

\def\dd{\partial}

\def\al{{\alpha}}
\def\be{{\beta}}
\def\Ga{{\Gamma}}
\def\ga{{\gamma}}
\def\de{{\delta}}

\def\la{{\lambda}}
\def\si{{\sigma}}
\def\om{{\omega}}

\def\hsi{{\hat\si}}
\def\arf{\si}

\def\st{\,\,\big|\,\,}
\def\<{\langle}
\def\>{\rangle}
\def\ie{i.e.\xspace}
\def\bs{\backslash}

\let\ge=\geqslant
\let\le=\leqslant

\DeclareMathOperator{\Aut}{Aut}
\DeclareMathOperator{\Arf}{Arf}

\DeclareMathOperator{\id}{id}
\DeclareMathOperator{\Hol}{Hol}
\DeclareMathOperator{\MathOpPSL}{PSL}
\DeclareMathOperator{\Mod}{Mod}
\DeclareMathOperator{\MathOpSL}{SL}
\DeclareMathOperator{\MathOpSp}{Sp}

\DeclareMathOperator{\trace}{trace}

\let\Im=\undefined \DeclareMathOperator{\Im}{Im}
 
\let\mod=\undefined \DeclareMathOperator{\mod}{mod}

\begin{document}

\author[Sergey Natanzon]{Sergey Natanzon}
\address{Moscow State University, Korp.~A, Leninske Gory, 11899 Moscow, Russia}
\address{Institute of Theoretical and Experimental Physics, Moscow, Russia}
\address{Independent University of Moscow, Bolshoi Vlasevsky Pereulok, 11 Moscow, Russia}
\email{natanzon@mccme.ru}
\author[Anna Pratoussevitch]{Anna Pratoussevitch}
\address{Department of Mathematical Sciences\\ University of Liverpool\\ Peach Street \\ Liverpool L69~7ZL}
\email{annap@liv.ac.uk}
%\address{Mathematisches Institut\\ Universit\"at Bonn\\ Beringstra{\ss}e~1 \\ 53115 Bonn}
%\email{anna@math.uni-bonn.de}

\title[Higher Arf Functions and Moduli Space of Higher Spin Surfaces]
{Higher Arf Functions and Topology of the Moduli Space of Higher Spin Riemann Surfaces}

%\title[Topological Invariants of Gorenstein Singularities]
%{Topological Invariants of Gorenstein Quasi-Homogeneous Singularities}

%\title{Higher Arf Functions and Gorenstein Singularities}
%\title{Higher Spin Surfaces and Arf Functions}

\begin{date}  {\today} \end{date}

\thanks{Research partially supported by grants INTAS 05-7805, RFBR-07-01-00593,
NSh-709.2008.1, NWO-RFBR 047.011.2004.026 (RFBR: 05-02-89000-HBO\_a),
and SFB 611 (DFG)}

\begin{abstract}

We describe all connected components of the space of pairs~$(P,s)$,
where $P$ is a hyperbolic Riemann surface with finitely generated
fundamental group and $s$ is an $m$-spin structure on~$P$. We prove
that any connected component is homeomorphic to a quotient of
${\mathbb R}^d$ by a discrete group.

\smallskip\noindent
Our method is based on a description of an $m$-spin structure by an $m$-Arf function,
that is a map $\sigma:\pi_1(P,p)\rightarrow{\mathbb Z}/m{\mathbb Z}$
with certain geometric properties. We prove that the set of all $m$-Arf functions
has a structure of an affine space associated
with $H_1(P,{\mathbb Z}/m{\mathbb Z})$.
We describe the orbits of $m$-Arf functions under the action of the group of homotopy classes of surface autohomeomorphisms.
Natural topological invariants of an orbit are the unordered set of va\-lues of the $m$-Arf functions on the punctures
and the unordered set of values on the $m$-Arf-function on the holes.
We prove
%(Theorem~\ref{thm-is-type})
that for~$g>1$ the space of $m$-Arf functions with prescribed genus and prescribed (unordered) sets of va\-lues on punctures and holes is either connected
% if it is not empty
%if $m$ is odd or if $m$ is even and the value of the $m$-Arf function on one of the punctures or holes is even.
% Otherwise this space
or has two connected components distinguished by the Arf invariant $\de\in\{0,1\}$.
%(see Definition~\ref{def-arf-inv})
(See the results for~$g=1$ later in the paper.)

\smallskip\noindent

%We describe the space of $m$-spin structures
%on a Riemann surface as a finite affine space
%of $({\mathbb Z}/m{\mathbb Z})$-valued functions
%on the fundamental group of the surface.
%We apply this description to prove that any connected component of the space
%of $m$-spin structures on compact Riemann surfaces with finite number of
%punctures and holes is homeomorphic to a quotient
% of the vector space ${\mathbb R}^d$ by a discrete group action.

\end{abstract}


% AG9809138: Primary 14H10, 32G15; Secondary 32G81, 81T40, 14M30.

% AG9908085: Primary 14H10, 32G15, 81T40; Secondary 14N, 14M.

% Q-Gorenstein papers of mine: Primary 32S25; Secondary 14J17, 14J60

% 32S25: Surface and hypersurface singularities

% 30F10: Compact Riemann surfaces and uniformization

% 14J17: Singularities

% 14J60: Vector bundles on surfaces and higher-dimensional varieties, and their
% moduli

% 14H10: Families, moduli (algebraic)

% 14H15: Families, moduli (analytic)

% 14M: Special varieties

% 14N: Projective and enumerative geometry

% 32G: Deformations of analytic structures

% 32G15: Moduli of Riemann surfaces, Teichmueller theory

% 32G81: Applications to physics

% 81T40: Two-dimensional field theories, conformal field theories, etc.


% singularities

% AG9809138, AG9908085: Algebraic curves, moduli, higher spin curves.

\maketitle

\section{Introduction}

% comment: autohomeo
% comment: MSC !

In this paper we study $m$-spin structures, \ie complex line bundles such that
the $m$-th tensor power is isomorphic to the cotangent bundle of the surface,
on hyperbolic Riemann surfaces with finitely generated fundamental group.

\myskip
The classical spin structures (theta characteristics) on compact Riemann surfaces
play an important role in algebraic geometry since Riemann~\cite{R}.
%(pp.~212, 487)
Their modern interpretation and classification as complex line
bundles such that the tensor square is isomorphic to the cotangent
bundle of the surface was given by Atiyah~\cite{A} and
Mumford~\cite{Mu}, a topological interpretation of their results was
given by Johnson~\cite{J}. They showed a connection between the set
of spin bundles on a surface~$P$ and the affine space of quadratic
(with respect to the index of intersection) forms
$H_1(P;\zd)\to\zd$. (Here $\zm=\z/m\z$.) Classification of classical
spin structures on non-compact Riemann surfaces and the
corresponding moduli space were studied in~\cite{N1989moduli,
Nbook}.

\myskip
We consider the moduli space of $m$-spin structures on Riemann surfaces,
that is the space of pairs~$(P,s)$,
where $P$ is a Riemann surface and $s$ is an $m$-spin structure on~$P$.
This moduli space plays an important role in mathematical physics~\cite{Wi,JKV:2001}
and singularity theory~\cite{Dolgachev:1983}.
%These moduli spaces were studied in \cite{JKV:2001, CCC:2004} and other papers.
%The number of connected components of the moduli space of $m$-spin structures on Riemann surfaces with punctures (but without holes) was computed in~\cite{Jarvis:2000}
%using methods of algebraic geometry.
The number of connected components of the moduli space of $m$-spin
structures on Riemann surfaces with punctures (but without holes)
was determined in~\cite{Jarvis:2000}, using an interpretation of the
quotients of $m$-spin structures on~$P$ as points of the Jacobian
of~$P$.

\myskip
%%The list of connected components of the space of compact $m$-spin Riemann surfaces
%%with punctures was found in~\cite{Jarvis:2000} using methods of algebraic geometry.
%%We give a description of the connected components of the moduli space
%%in terms of $m$-Arf functions and find the number of connected components (Theorem~\ref{moral}).
%%For surfaces without holes this number was computed in~\cite{Jarvis:2000}
%%using methods of algebraic geometry.
%%[Our methods are different.
%%Our methods are topological.
%%With our methods
%We use topological methods to study the moduli spaces of $m$-spin structures
%on hyperbolic Riemann surfaces with punctures and holes.
%We are able to determine the connected components of the moduli space.
%Moreover, with our methods we can describe the topology of the connected components.
%%We identify them with...]
%%In our paper we give an independent topological description of the connected
%%components using the topological invariants of the higher Arf functions.
Our method uses Fuchsian groups and allows us to work with the fundamental group of the surface rather than with the homology group.
We assign to any $m$-spin structure a set of topological invariants that describes the behavior of the $m$-spin structure on the basis of the fundamental group.
We then show that this set of topological invariants determines a connected component of the moduli space of $m$-spin structures.
Moreover we prove
%%(Theorem~\ref{top-type-comp})
that any connected component is homeomorphic to the space of the form $\r^d/\Mod$,
where $\Mod$ is a discrete group acting on~$\r^d$.

\myskip
The main technical tool is the following:
We assign (Theorem~\ref{thm-corresp}) to any $m$-spin structure on a surface~$P$
a unique function on the space of homotopy classes of simple contours on~$P$
with values in $\zm$, the associated $m$-Arf function.
In contrast with the case of classical spin structures on compact Riemann surfaces
these functions are not commutative and do not induce quadratic forms on the homology group.
For this reason some of the proofs are more computational and less geometrical than Johnson's results in~\cite{J}.
% (we call them higher Arf functions)

\myskip\noindent
The $m$-Arf functions are described by simple geometric properties:

\myskip\noindent
{\bf Definition:}
We denote by $\pi_1^0(P,p)$ the set of all non-trivial elements of $\pi_1(P,p)$
that can be represented by simple contours.
%that either do not belong to the kernel of the intersection form or are homologous to a hole or a puncture.
An {\it $m$-Arf function} is a function
$$\arf:\pi_1^0(P,p)\to\zm$$
satisfying the following conditions
\begin{enumerate}[1.]
\item
$\arf(bab^{-1})=\arf(a)$ for any elements~$a,b\in\pi_1^0(P,p)$,
% such that the element~$bab^{-1}$ is in~$\pi_1^0(P,p)$,
\item
$\arf(a^{-1})=-\arf(a)$ for any element~$a\in\pi_1^0(P,p)$,
\item
$\arf(ab)=\arf(a)+\arf(b)$
for any
% elements~$a,b\in\pi_1^0(P,p)$ such that the element~$ab$ is in~$\pi_1^0(P,p)$ and the
elements~$a$ and~$b$ which can be represented by a pair of simple contours in $P$
intersecting in exactly one point~$p$ with $\<a,b\>\ne0$,
\item
$\arf(ab)=\arf(a)+\arf(b)-1$ for any
elements~$a,b\in\pi_1^0(P,p)$ such that the element~$ab$ is
in~$\pi_1^0(P,p)$ and the elements~$a$ and~$b$ can be
represented by a pair of simple contours in $P$ intersecting in
exactly one point~$p$ with $\<a,b\>=0$ and placed in a
neighborhood of the point~$p$ as shown in
Figure~\ref{fig-neg-pair}.
%\ie in such a way that the oriented contours $a$, $b$, and
%$(ab)^{-1}$ are freely homotopic to pairwise non-intersecting simple contours
%with orientation opposite to the one induced by the complex structure of the sphere with three
%holes that they cut out of $P$.
\end{enumerate}

%\myskip\noindent
%The simple contours are primitive elements of~$\pi_1(P,p)$,
%thus we can continue~$\arf$ to~$\pi_1(P,p)$.
% by~$\arf(ka)=k\arf(a)$ retaining the properties~1.--3.

\myskip\noindent
In order to formulate our main results we need to give some definitions and notation.
We say that a hyperbolic Riemann surface is of type $(g,l_h,l_p)$
if the surface is obtained from a surface of genus~$g$
by removing $l_h$~disks and $l_p$~points.
%We say that a hyperbolic Riemann surface of type $(g,l_h,l_p)$
%is a surface of genus~$g$
%with $l_h$ (hyperbolic) holes and $l_p$ (parabolic) punctures.

\myskip\noindent
{\bf Definition:}
Let $P$ be a hyperbolic Riemann surface of type $(g,l_h,l_p)$.
Let $\arf:\pi_1^0(P,p)\to\zm$ be an $m$-Arf function on~$P$.
For~$g>1$ and even~$m$ the Arf invariant $\de=\de(P,\arf)$ of the $m$-Arf function~$\arf$
equals zero if there is a standard basis
$$\{a_i,b_i~(i=1,\dots,g),c_i~(i=g+1,\dots,n)\}$$
of the fundamental group $\pi_1(P,p)$
(compare to Definition~\ref{def-standard-basis}) such that
$$\sum\limits_{i=1}^g(1-\arf(a_i))(1-\arf(b_i))=0~\mod~2$$
and equals one otherwise.
For~$g>1$ and odd $m$ we set $\de=0$.
(For~$g=1$, see the definition of the Arf invariant later in the paper.)
% The Arf invariant is well-defined.
The type of the $m$-Arf function~$\arf$ is the tuple
$$(g,\de,\num^h_0,\dots,\num^h_{m-1},\num^p_0,\dots,\num^p_{m-1}),$$
where $\de$ is the Arf invariant of~$\arf$,
and $\num^h_j$ resp.\ $\num^p_j$ for $j\in\zm$
is the number of holes resp.\ punctures,
such that the Arf function~$\arf$ attains the value~$j$
on the corresponding element of the standard basis.

\myskip\noindent
{\bf Definition:}
We denote by
$$\Smgll(t)=\Smgll(g,\de,\num^h_0,\dots,\num^h_{m-1},\num^p_0,\dots,\num^p_{m-1})$$
the set of all $m$-spin structures on Riemann surfaces
of type $(g,l_h,l_p)$ such that the associated $m$-Arf function is of type
$t=(g,\de,\num^h_0,\dots,\num^h_{m-1},\num^p_0,\dots,\num^p_{m-1})$.
Here
$$l_h=\sum_{j\in\zm}\,\num^h_j\quad\text{and}\quad l_p=\sum_{j\in\zm}\,\num^p_j.$$

\myskip\noindent
The following Theorem summarizes the main results:

\begin{thm}

\noindent
\begin{enumerate}[1)]
\item
Two $m$-spin structures are in the same connected component
of the space of all $m$-spin structures on hyperbolic Riemann surfaces
iff they are of the same type.
In other words, the connected components of the space of all $m$-spin
structures are those sets $\Smgll(t)$
that are not empty.
\item
The set $\Smgll(t)$ is not empty iff
$t=(g,\de,\num^h_0,\dots,\num^h_{m-1},\num^p_0,\dots,\num^p_{m-1})$
has the following properties:
\begin{enumerate}[(a)]
\item
If $g>1$ and $m$ is odd, then $\de=0$.
\item
If $g>1$ and $m$ is even and $\num^h_j+\num^p_j\ne0$ for some even~$j\in\zm$, then $\de=0$.
\item
If $g=1$ then $\de$ is a divisor of~$m$ and~$\gcd(\{j+1\st\num^h_j+\num^p_j\ne0\})$.
\item
The following degree condition is satisfied
$$\sum\limits_{j\in\zm}\,j\cdot(\num^h_j+\num^p_j)=(2-2g)-(l_h+l_p).$$
\end{enumerate}
\item
Any connected component~$\Smgll(t)$ of the space of all $m$-spin structures on hyperbolic Riemann surfaces of type~$(g,l_h,l_p)$
is homeomorphic to a quotient of the space~$\r^{6g+3l_h+2l_p-6}$
by a discrete action of a certain subgroup of the modular group (see subsection~\ref{components} for details).
\end{enumerate}
\end{thm}

\myskip
Let us outline the main ideas and methods of the proof.
We use the higher Arf functions to describe the topology of the
moduli space of $m$-spin bundles.

\myskip
We first outline the construction assigning an $m$-Arf function
to an $m$-spin structure on $P=\hyp/\lat$,
where $\hyp$ is the hyperbolic plane and $\lat$ is a Fuchsian group without elliptic elements.
The construction is based on the topological properties
of the group $\PSL\cong\Aut(\hyp)$.
This group has a unique connected $m$-fold co\-ve\-ring $\Gm\to\PSL$.
Thereby there is a 1-1-correspondence between $m$-spin structures
on $P=\hyp/\lat$ and lifts into $\Gm$ of the group $\lat$,
\ie subgroups $\lats$ of $\Gm$ such that the restriction
of the covering map $\Gm\to\PSL$ to $\lats$ is an isomorphism $\lats\to\lat$~\cite{Mi}.
We prove that the preimage in $\Gm$ of the set of all hyperbolic and parabolic
elements of $\PSL$ has $m$ connected components,
which we identify with elements of the group $\zm$.
This correspondence induces a map $\arf:\pi_1^0(P,p)\to\zm$.
The geometric properties of this map follow from the discreteness
criterion~\cite{Nbook} for subgroups of $\PSL$.

\myskip
We prove that the the set of all such functions has a structure of an affine space
associated with $H^1(P;\zm)$.
We describe the orbits of $m$-Arf functions under the action of the group of homotopy classes of surface autohomeomorphisms.
Natural topological invariants of an orbit are the unordered sets of va\-lues of the $m$-Arf functions on the punctures resp.\ holes.
We prove (in the case~$g>1$)
%(Theorem~\ref{thm-is-type})
that the space of $m$-Arf functions with prescribed genus
%of the surface
and prescribed (unordered) sets of va\-lues on punctures resp. holes
is either connected
% if it is not empty
%if $m$ is odd or if $m$ is even and the value of the $m$-Arf function
% on one of the punctures or holes is even.
% Otherwise this space
or has two connected components distinguished by the Arf invariant $\de\in\{0,1\}$.
%(see Definition~\ref{def-arf-inv})

\myskip
The paper is organized as follows:
In section~\ref{sec-PSL} we study the covering groups $\Gm$ of the group
$\PSL$, and in particular the algebraic properties of the preimages in $\Gm$
of hyperbolic and parabolic elements of $\PSL$.
In section~\ref{sec-spins-and-lifts} we explore the connection
between $m$-spin structures on a Riemann surface~$P=\hyp/\lat$
and lifts into the covering $\Gm$ of the group $\lat$.
We assign to any lift a function induced by a decomposition of the covering
$\Gm$ into sheets and choosing a numeration of the sheets and study properties
of these functions.
In section~\ref{sec-m-arf} we define $m$-Arf functions.
We prove that there is a 1-1-correspondence between the set of $m$-Arf functions on~$P=\hyp/\lat$
and the set of functions associated to the lifts of $\lat$ via the numeration of the covering sheets.
Hence these two sets are also in 1-1-correspondence with the set of $m$-spin structures on~$P$.
Moreover we show in this section using the explicit description of the Dehn generators
of the group of homotopy classes of surface autohomeomorphisms
that the set of all $m$-Arf functions on a surface~$P$ has a structure of an affine space.
In the last section we find topological invariants of $m$-Arf functions
and prove that they describe the connected components of the moduli space.
Furthermore we show using a version of Theorem of Fricke and Klein
that any connected component is homeomorphic to the space of the form $\r^d/\Mod$,
where $\Mod$ is a discrete group acting on~$\r^d$.

\myskip
Part of this work was done during the stays at Max-Planck-Institute in Bonn and at IHES.
We are grateful to the both institutions for their hospitality and support.
We would like to thank E.B.~Vinberg for many useful discussions related to this work.
%We would like to thank Walter Neumann for useful conversations related to this work.
We would like to thank the referees for their valuable remarks and suggestions.

\section{The group ${\rm PSL}(2,\mathbb R)$ and its coverings}

\label{sec-PSL}

We consider the universal cover $\tG=\tPSL$ of the Lie group
$$G=\PSL=\SL/\{\pm1\},$$
the group of orien\-ta\-tion-preserving isometries of the hyperbolic plane.
Here our model of the hyperbolic plane is the upper half-plane $\hyp=\{z\in\c\st\Im(z)>0\}$
and the action of an element $[\abcd]\in\PSL$ on~$\hyp$ is by
$$z\mapsto\frac{az+b}{cz+d}.$$
Here we denote by $[A]=[\abcd]\in\PSL$ the equivalence class of a matrix $A=\abcd\in\SL$.

\begin{notation}
$\s^1=\{z\in\c\st|z|=1\}\subset\c$, $\r_+=\{x\in\r\st x>0\}$, $\zm=\z/m\z$.
\end{notation}

\subsection{Description of elements in $G={\rm PSL}(2,\mathbb R)$}

Elements of $\PSL$ can be classified with respect to the fixed point behavior
of their action on~$\hyp$.
An element is called {\it hyperbolic} if it has two fixed points,
which lie on the boundary $\dd\hyp=\r\cup\{\infty\}$ of~$\hyp$.
A hyperbolic element with fixed points~$\al$, $\beta$ in~$\r$ is of the form
$$
  \tau_{\al,\be}(\la)
  :=\left[\frac{1}{(\al-\be)\cdot\sqrt{\la}}\cdot\begin{pmatrix}\la\al-\be&-(\la-1)\al\be\\ \la-1&\al-\la\be\end{pmatrix}\right],
$$
where~$\la>0$.
A hyperbolic element with one fixed point in~$\infty$ is of the form
% comment: \lim\limits_{\al\to\infty}\tau_{\al,\be}
$$
  \tau_{\infty,\be}(\la)
  :=\left[\frac{1}{\sqrt{\la}}\cdot\begin{pmatrix}\la&-(\la-1)\be\\ 0&1\end{pmatrix}\right]
$$
or
% comment: \lim\limits_{\be\to\infty}\tau_{\al,\be}
$$
  \tau_{\al,\infty}(\la)
  :=\left[\frac{1}{\sqrt{\la}}\cdot\begin{pmatrix}1&(\la-1)\al\\ 0&\la\end{pmatrix}\right],
$$
where $\al$ resp.\ $\be$ is the real fixed points and $\la>0$.
The parameter $\la>0$ is called the {\it shift parameter}.
The {\it axis} $\ell(g)$ of the element $g=\tau_{\al,\be}(\la)$ is the geodesic between the fixed points $\al$ and $\be$,
oriented from $\be$ to $\al$ if $\la>1$ and from $\al$ to $\be$ if $\la<1$.
The element $g=\tau_{\al,\be}(\la)$ preserves the geodesic $\ell(g)$ and moves
the points on this geodesic in the direction of the orientation.
We call a hyperbolic element $\tau_{\al,\be}(\la)$ with $\la>1$ {\it positive} if $\al<\be$.
The map $\la\mapsto\tau_{\al,\be}(\la)$ defines a homomorphism $\r_+\to G$
(with respect to the multiplicative structure on $\r_+$).
We have
$$(\tau_{\al,\be}(\la))^{-1}=\tau_{\al,\be}(\la^{-1})=\tau_{\be,\al}(\la).$$

\myskip
An element is called {\it parabolic} if it has one fixed point, which is on the boundary~$\dd\hyp$.
A parabolic element with real fixed point $\al$ is of the form
$$
  \pi_{\al}(\la)
  :=\left[\begin{pmatrix}1-\la\al&\la\al^2\\ -\la&1+\la\al\end{pmatrix}\right].
$$
A parabolic element with fixed point $\infty$ is of the form
$$
  \pi_{\infty}(\la)
  :=\left[\begin{pmatrix}1&\la\\ 0&1\end{pmatrix}\right].
$$
We call a parabolic element $\pi_{\al}(\la)$ {\it positive} if $\la>0$.
The map $\la\mapsto\pi_{\al}(\la)$ defines a homomorphism $\r\to G$.
(with respect to the additive structure on $\r$).
We have
$$(\pi_{\al}(\la))^{-1}=\pi_{\al}(-\la).$$

\myskip
An element that is neither hyperbolic nor parabolic is called {\it elliptic}.
It has one fixed point that is in $\hyp$.
Given a base-point $x\in\hyp$ and a real number $\varphi$,
let $\rho_x(\varphi)\in G$ denote the rotation through angle $\varphi$ about the point $x$.
Any elliptic element is of the form $\rho_x(\varphi)$, where $x$ is the fixed point.
Thus we obtain a $2\pi$-periodic homomorphism $\rho_x:\r\to G$
(with respect to the additive structure on $\r$).
We have
$$\rho_x(\varphi+2\pi)=\rho_x(\varphi)\quad\hbox{and}\quad(\rho_x(\varphi))^{-1}=\rho_x(-\varphi).$$
For the fixed point $x=i\in\hyp$ we have
$$
  \rho_i(\varphi)
  =\left[\begin{pmatrix} \cos\frac{\varphi}2&-\sin\frac{\varphi}2\\ \sin\frac{\varphi}2&\cos\frac{\varphi}2\end{pmatrix}\right].
$$
For a fixed point $x\in\hyp\bs\{i\}$ we obtain $\rho_x(\varphi)=\tau\circ\rho_i(\varphi)\circ\tau^{-1}$,
where $\tau$ is the hyperbolic element in~$G$ such that $\tau(i)=x$.

\subsection{Coverings $\Gm$ of $G$}

As topological space $\PSL$ is homeomorphic to the open solid torus $\s^1\times\c$.
A homeomorphism $\PSL\to\s^1\times\c/\{\pm1\}$ can be given explicitely, see subsection~\ref{subsec-sheets}.
%A homeomorphism $\PSL\to\s^1\times\c/\{\pm1\}$ can be given explicitely by
%$$H:\left[\ABCD\right]\mapsto\left(\left(\frac{(a+d)+i(b-c)}{|(a+d)+i(b-c)|}\right)^2,\left[\frac{(a-d)+i(b+c)}{2}\right]\right).$$

\myskip
The fundamental group of the open solid torus $G$ is infinite cyclic.
Therefore, for each natural number $m$ there is a unique connected $m$-fold covering
$$\Gm=\tG/(m\cdot Z(\tG))$$
of $G$, where $\tG$ is the universal covering of~$G$ and $Z(\tG)$ is the centre of~$\tG$.
For $m=2$ this is the group $G_2=\SL$.

%\myskip
%The map $\mu:G\to\s^1$ defined as the composition of the homeomorphism~$H$
%and the projection onto the factor $\s^1$
%maps an element $[A]=[\abcd]\in G$ to a unit complex number~$\mu([A])=e^{i\psi}$ with
%$$\tan\frac{\psi}{2}=\frac{b-c}{a+d}.$$
%We shall refer to the number $\psi$ as the argument of the element~$[A]$.
%The map $\mu:G\to\s^1$ lifts to the unique map $\varphi:\tG\to\r$ of the universal covers
%such that the following diagram commutes
%$$
%  \begin{CD}
%   \tG            @>\varphi>> \r         \\
%   @V{}VV         @VV{}V \\
%   G            @>\mu>> \s^1          \\
%  \end{CD}
%$$
%and $\varphi(\tilde e)=0$, where $\tilde e$ is the identity element in~$\tG$.
%Here the map $\r\to\s^1$ is the universal covering map given by $x\mapsto e^{ix}$.

%\begin{mydef}
%We define level functions~$\lev:\tG\to\z$ on~$\tG$
%and~$\levm:\Gm\to\zm$ on $\Gm=\tG/(m\cdot Z(\tG))$ by
%$$
%  \lev(g)=k
%  \quad\text{for}\quad
%  g\in\tG
%  \quad\text{if}\quad
%  \varphi(g)\in(-\pi+2\pi k,\pi+2\pi k]
%$$
%and
%$$\levm(g~\mod~(m\cdot Z(\tG)))=\lev(g)~\mod~m\quad\text{for}\quad g\in\tG.$$
%We recall that $\zm=\z/m\z$.
% All equations involving~$\levm$ are to be understood as equations in~$\zm$, \ie equations modulo~$m$.
%\end{mydef}

\myskip
Here is another description of the covering groups $\Gm$ of $G=\PSL$, which fixes a group structure.
Let $\Hol(\hyp,\c^*)$ be the set of all holomorphic functions $\hyp\to\c^*$.

\begin{prop}
\label{fract-autforms}
The $m$-fold covering group $\Gm$ of~$G$ can be described as
$$\{(g,\de)\in G\times\Hol(\hyp,\c^*)\st \de^m(z)=g'(z)~{\rm for~all}~z\in\hyp\}$$
with multiplication
$$(g_2,\de_2)\cdot(g_1,\de_1)=(g_2\cdot g_1,(\de_2\circ g_1)\cdot\de_1).$$
\end{prop}

\begin{proof}
Let $X$ be the subspace of $G\times\Hol(\hyp,\c^*)$ in question.
One can check that the space $X$ is connected
and that the map $X\to G$ given by $(\ga,\de)\mapsto\ga$ is an $m$-fold covering of~$G$.
Hence the coverings $X\to G$ and $\Gm\to G$ are isomorphic.
One can check that the operation described above defines a group structure on~$X$
and that the covering map $X\to G$ is a homomorphism with respect to this group structure.
\end{proof}

\begin{rem}
This description of $\Gm$ is inspired by the notion of automorphic
differential forms of fractional degree, introduced by J.~Milnor in~\cite{Mi}.
For a more detailed discussion of this fact see~\cite{Lion:Vergne}, section~1.8.
\end{rem}

Elements of $\tG$ resp.~$\Gm$ can also be classified with respect to the fixed point
behavior of action on $\hyp$ of their image in~$G=\PSL$.

\begin{mydef}
We say that an element of~$\tG$ resp.~$\Gm$ is {\it hyperbolic,} {\it parabolic,} resp.\ {\it elliptic} if its image in~$G$ has this property.
\end{mydef}

%We also define the absolute value of the trace for an element of~$\Gm$
%as the absolute value of the trace of its image in~$G$.

\subsection{One-Parameter-Subgroups of $G$ and $\tG$}

\myskip
The homomorphisms
$$\tau_{\al,\be}:\r_+\to G,\quad \pi_{\al}:\r\to G,\quad\text{resp.}\quad \rho_x:\r\to G$$
define one-parameter-subgroups in the group~$G=\PSL$.

\myskip
Each of the homomorphisms
$\tau_{\al,\be}:\r_+\to G$, $\pi_{\al}:\r\to G$, resp.\ $\rho_x:\r\to G$
lifts to a unique homomorphism
$$t_{\al,\be}:\r_+\to\tG,\quad p_{\al}:\r\to\tG\quad\text{resp.}\quad r_x:\r\to\tG$$
into the universal covering group.
The elements $t_{\al,\be}(\la)$, $p_{\al}(\la)$, resp.\ $r_x(\xi)$
are hyperbolic, parabolic, resp.\ elliptic.

%\myskip
%For hyperbolic elements we obtain for all $\la\in\r_+$ that
%$$\varphi(t_{\al,\beta}(\la))\in(-\pi,\pi)$$
%because of $\mu(\tau_{\al,\beta}(1))=1$ and $\mu(\tau_{\al,\beta}(\la))\ne-1$.

%\myskip
%Similarly, for parabolic elements we obtain for all $\la\in\r$ that
%$$\varphi(p_{\al}(\la))\in(-\pi,\pi)$$
%because of $\mu(\pi_{\al}(0))=1$ and $\mu(\pi_{\al}(\la))\ne-1$.

%\myskip
%Now let us consider elliptic elements.
%For the trace of the rotation $\rho_x(\xi)$ we have
%$$|\trace(\rho_x(\xi))|=|\trace(\rho_i(\xi))|=|2\cos(\xi/2)|,$$
%hence $\mu(\rho_x(\xi))=-1$ if and only if $\xi\in\pi+2\pi\cdot\z$.
%Because of $\mu(\rho_x(0))=1$ this implies for $\xi\in(-\pi+2\pi k,\pi+2\pi k)$ that
%$$\varphi(r_x(\xi))\in(-\pi+2\pi k,\pi+2\pi k).$$
%We obtain $\lev(r_x(\xi))=k$ for $\xi\in(-\pi+2\pi k,\pi+2\pi k]$.

\subsection{The Centre of the Group $\tG$}

\myskip
Since $\rho_x(2\pi\ell)=\id$ for~$\ell\in\z$,
it follows that the lifted element $r_x(2\pi\ell)$ belongs to the centre $Z(\tG)$ of~$\tG$.
Note that this element $r_x(2\pi\ell)$ depends continuously on $x$.
But the centre of~$\tG$ is discrete, so this element must remain constant,
thus~$r_x(2\pi\ell)$ does not depend on~$x$.
The centre~$Z(\tG)$ of~$\tG$ is equal to the pre-image of the identity element under the projection~$\tG\to G$, hence
$$Z(\tG)=\{r_x(2\pi\ell)\st\ell\in\z\}.$$
Let
$$u=r_x(2\pi)$$
for some (and hence for any)~$x$ in~$\hyp$.
The element~$u$ is one of the two generators of the centre of~$\tG$,
since any other element of the centre is of the form~$r_x(2\pi\ell)=(r_x(2\pi))^{\ell}$.

%\myskip
%We would also like to point out that for the lift of an elliptic element $\rho_x(2\pi/p)$ of finite order~$p$ we have
%$$(r_x(2\pi/p))^p=r_x(2\pi)=u.$$

\subsection{Decomposition of the Subset of Hyperbolic and Parabolic Elements in $\tG$ and $\Gm$ into Sheets}
\label{subsec-sheets}

\myskip
Let $\Xi$ be the subset of~$G=\PSL$ that consists of all hyperbolic and parabolic elements of~$G$ (including the identity element).
The space~$G$ is homeomorphic to the open solid torus~$\s^1\times\c$.
In~\cite{JN} Jankins and Neumann give an explicite homeomorphism (see~\cite{JN}, Apendix)
and describe the image of the subset~$\Xi$ under this homeomorphism (see~\cite{JN}, \S~1, Figure~1).
From this description it follows in particular that the subset~$\Xi$ is simply connected.
The pre-image~$\tilde\Xi$ of the subset~$\Xi$ in~$\tG$ consists of infinitely many connected components.
Each connected component of the subset~$\tilde\Xi$ contains one and only one pre-image of the identity element of~$G$,
\ie one and only one element of the centre of~$\tG$.
The elements of the one-parameter-subgroups~$t_{\al,\be}(\la)$ resp.~$p_{\al}(\la)$ are contained in the same connected component of the subset~$\tilde\Xi$
as the identity element~$\tilde e$ of the group~$\tG$.

%The homeomorphism $\PSL\to\s^1\times\c/\{\pm1\}$ can be given explicitely as
%$$H:\left[\ABCD\right]\mapsto\left(\left(\frac{(a+d)+i(b-c)}{|(a+d)+i(b-c)|}\right)^2,\left[\frac{(a-d)+i(b+c)}{2}\right]\right).$$
%Here we denote by $[A]=[\abcd]\in\PSL$ the equivalence class of a matrix $A=\abcd\in\SL$, and by $[z]\in\c/\{\pm1\}$ the equivalence class of~$z\in\c$.
%Under the homeomorphism~$H:G\to\s^1\times\c/\{\pm1\}$ the identity element~${1\,0\choose 0\,1}$ corresponds to the point~$(1,0)\in\s^1\times\c$.

\begin{mydef}
If a hyperbolic or parabolic element of~$\tG$ is contained in the same connected component as the central element~$u^k$, $k\in\z$,
we say that the element is {\it at the level}~$k$ and set the {\it level function}~$\lev$ on this element to be equal to~$k$.
Any hyperbolic resp.\ parabolic
%resp.\ elliptic
element in~$\tG$ is of the form
$t_{\al,\be}(\la)\cdot u^k$ resp.\ $p_{\al}(\la)\cdot u^k$.
% resp.\ $r_x(\xi)$.
For elements written in this form we have
$$\lev(t_{\al,\be}(\la)\cdot u^k)=k,\quad\lev(p_{\al}(\la)\cdot u^k)=k.$$
%and $\lev(r_x(\xi))=k$ iff $\xi\in(-\pi+2\pi k,\pi+2\pi k]$.
\end{mydef}

\begin{mydef}
We define the {\it level function}~$\levm$ on the hyperbolic and parabolic elements of~$\Gm=\tG/(m\cdot Z(\tG))$ by
$$\levm(g~\mod~(m\cdot Z(\tG)))=\lev(g)~\mod~m\quad\text{for}\quad g\in\tG.$$
We recall that $\zm=\z/m\z$.
All equations involving~$\levm$ are to be understood as equations in~$\zm$, \ie equations modulo~$m$.
\end{mydef}

\subsection{Properties of the multiplication of hyperbolic and parabolic elements in $\Gm$}

\label{mult-tG}

In this subsection we first study (Lemma~\ref{lem-inv} and~\ref{lem-conj})
the behavior of the level functions $\levm$ under inversion and conjugation.
The main results of this subsection (Lemma~\ref{lem-product-of-hyps}, \ref{lem-product-of-hyp-and-par}, and \ref{lem-product-of-pars})
are statements about the behavior of $\levm$ under multiplication.

\myskip
In this subsection let us denote by $[\cdot]$ the image of an element in~$\tG$ under the covering map $\tG\to\Gm$.

\begin{lem}
\label{lem-inv}
The equation $\levm(A^{-1})=-\levm(A)$ is satisfied for any hyperbolic or parabolic element~$A$ in~$\Gm$.
\end{lem}

\begin{proof}
The hyperbolic resp.\ parabolic element $A$ is of the form $[t_{\al,\be}(\la)\cdot u^k]$ resp.\ $[p_{\al}(\la)\cdot u^k]$.
If $A=[t_{\al,\be}(\la)\cdot u^k]$
then $A^{-1}=[t_{\al,\be}(\la^{-1})\cdot u^{-k}]$
and $\levm(A^{-1})=-k=-\levm(A)$.
If $A=[p_{\al}(\la)\cdot u^k]$
then $A^{-1}=[p_{\al}(-\la)\cdot u^{-k}]$
and $\levm(A^{-1})=-k=-\levm(A)$.
\end{proof}

\begin{lem}
\label{lem-conj}
For any hyperbolic or parabolic element~$A$ and any element~$B$ in~$\Gm$ we have
$$\levm(B\cdot A\cdot B^{-1})=\levm(A).$$
\end{lem}

\begin{proof}
%The expression $\levm(BAB^{-1})-\levm(A)$ is invariant under multiplication of $A$ or~$B$ with the central element~$[u]$,
%hence it is sufficient to prove the statement for elements $A$ and~$B$ with $\levm(A)=\levm(B)=0$.
%The element $A$ in $\Gm$ with $\levm(A)=0$ is of the form $[t_{\al,\be}(\la)]$, $[p_{\al}(\la)]$ resp.\ $[r_x(\varphi)]$ with $\varphi\in(-\pi,\pi]$.
%The element $A$ can be connected to the unit element in $\Gm$ via a path $\ga:I\to\Gm$, where $I$ is some closed interval, such that $\levm(\ga(t))=0$ for all $t\in I$.
The element~$B$ can be connected to the unit element in $\Gm$ via a path $\be:I\to\Gm$, where $I$ is some closed interval.
%For example, we can take the path $[t_{\al,\be}(t)]$, $t\in[1,\la]$, $[p_{\al}(t)]$, $t\in[0,\la]$, resp.\ $[r_x(t)]$, $t\in[0,\varphi]$.
%The path $\ga^B:I\to\Gm$ given by $\ga^B(t)=B\cdot\ga(t)\cdot B^{-1}$ connects the element $B\cdot A\cdot B^{-1}$ with the unit element.
The path $\ga:I\to\Gm$ given by
$$\ga(t)=\be(t)\cdot A\cdot(\be(t))^{-1}$$
connects the elements~$A$ and~$B\cdot A\cdot B^{-1}$.
%We have $|\trace(\ga^B(t))|=|\trace(B\cdot\gamma(t)\cdot B^{-1})|=|\trace(\gamma(t))|\ne0$ for all $t\in I$ and hence $\levm(\ga^B(t))=0$ for all $t\in I$,
%in particular $\levm(B\cdot A\cdot B^{-1})=0$.
We have
$$|\trace(\ga(t))|=|\trace(\be(t)\cdot A\cdot(\be(t))^{-1})|=|\trace(A)|\ne0$$
for all $t\in I$
and hence $\levm$ is constant along~$\ga$, in particular $\levm(B\cdot A\cdot B^{-1})=\levm(A)$.
(Here $|\trace X|$ for an element $X$ in $\Gm$ is defined as $|\trace\bX|$,
where $\bX$ is the projection of $X$ in $G$.)
\end{proof}

\begin{lem}
\label{lem-product-of-hyps}
For preimages $A$ resp.~$B$ in~$\Gm$
of the elements $\bA=\tau_{\infty,0}(\la_1)$ resp.\ $\bB=\tau_{\al,\be}(\la_2)$
for some $\la_1,\la_2>1$ and $\al,\be\in\r\bs\{0\}$
the difference
$$\levm(A\cdot B)-\levm(A)-\levm(B)$$
is equal to
\begin{enumerate}[$\bullet$]
\item $+1$ if
$$0<\frac{\la_1+\la_2}{1+\la_1\la_2}\cdot\be<\al<\be,$$
\item $-1$ if
$$\be<\al\le\frac{\la_1+\la_2}{1+\la_1\la_2}\cdot\be<0,$$
\item $0$ otherwise.
\end{enumerate}
\end{lem}

\myskip
One of the consequences of Lemma~\ref{lem-product-of-hyps} is the following corollary:

\begin{cor}
\label{cor-product-of-crossing-hyps}
If the axes of two hyperbolic elements $A$ and~$B$ in~$\Gm$ intersect,
%or the elements $A$ and $B$ are oriented differently (\ie $A$ is positive while $B$ is negative or $B$ is positive while $A$ is negative),
then
$$\levm(A\cdot B)=\levm(A)+\levm(B).$$
\end{cor}

\myskip
This corollary can be shown in a more geometrical way:

\myskip
Let $\ell_A$ resp.~$\ell_B$ be the axes of~$A$ resp.~$B$.
Let $x$ be the intersection point of~$\ell_A$ and~$\ell_B$.
Any hyperbolic transformation with the axis~$\ell_A$ is a product of a rotation by~$\pi$ at some point~$y\ne x$ on~$\ell_A$ and a rotation by~$\pi$ at the point~$x$.
Similarly any hyperbolic transformation with the axis~$\ell_B$ is a product of a rotation by~$\pi$ at the point~$x$ and a rotation by~$\pi$ at some point~$z\ne x$ on~$\ell_B$.
Hence the product of any hyperbolic transformation with the axis~$\ell_A$ and any hyperbolic transformation with the axis~$\ell_B$
is a product of a rotation by~$\pi$ at a point~$y\ne x$ on~$\ell_A$ and a rotation by~$\pi$ at a point~$z\ne x$ on~$\ell_B$,
\ie it is a hyperbolic transformation with an axis going through the points~$y$ and~$z$.
Thus the product of two hyperbolic elements with distinct but intersecting axes is always a hyperbolic element.

\myskip
Assume without loss of generality that the elements~$A$ and~$B$ satisfy the conditions~$\levm(A)=\levm(B)=0$.
We want to show that~$\levm(AB)=0$.
Let us deform the elements~$A$ and~$B$ by decreasing their shift parameters,
then their product tends to the identity element but remains hyperbolic, hence $\levm(AB)$ is equal to the value of~$\levm$ at the identity, \ie $\levm(AB)=0$.

\myskip
However, we can not give a similarly illuminating geometrical proof for the other parts of Lemma~\ref{lem-product-of-hyps}
and will prove these cases by computation.
To this end we will need the following definition:

\myskip
The homeomorphism $\PSL\to\s^1\times\c/\{\pm1\}$ can be given explicitely as
$$H:\left[\ABCD\right]\mapsto\left(\left(\frac{(a+d)+i(b-c)}{|(a+d)+i(b-c)|}\right)^2,\left[\frac{(a-d)+i(b+c)}{2}\right]\right).$$
Here we denote by $[A]=[\abcd]\in\PSL$ the equivalence class of a matrix $A=\abcd\in\SL$, and by $[z]\in\c/\{\pm1\}$ the equivalence class of~$z\in\c$.
The map $\mu:G\to\s^1$ defined as the composition of the homeomorphism~$H$
and the projection onto the factor $\s^1$
maps an element $[A]=[\abcd]\in G$ to a unit complex number~$\mu([A])=e^{i\psi}$ with
$$\tan\frac{\psi}{2}=\frac{b-c}{a+d}.$$
We shall refer to the number $\psi$ as the {\it argument} of the element~$[A]$.
The map $\mu:G\to\s^1$ lifts to the unique map $\varphi:\tG\to\r$ of the universal covers
such that the following diagram commutes
$$
  \begin{CD}
   \tG            @>\varphi>> \r         \\
   @V{}VV         @VV{}V \\
   G            @>\mu>> \s^1          \\
  \end{CD}
$$
and $\varphi(\tilde e)=0$, where $\tilde e$ is the identity element in~$\tG$.
Here the map $\r\to\s^1$ is the universal covering map given by $x\mapsto e^{ix}$.
The map~$\varphi$ will help us to compute the level function.
For hyperbolic and parabolic elements we have
$$
  \lev(g)=k
  \quad\text{for}\quad
  g\in\tG
  \quad\text{if and only if}\quad
  \varphi(g)\in(-\pi+2\pi k,\pi+2\pi k).
$$

\begin{proof}
The expression $\levm(A\cdot B)-\levm(A)-\levm(B)$ is invariant under multiplication of $A$ or~$B$ with the central element~$[u]$,
hence it is sufficient to prove the statement for elements $A$ and~$B$ with $\levm(A)=\levm(B)=0$,
\ie for $A=[t_{\infty,0}(\la_1)]$ and $B=[t_{\al,\be}(\la_2)]$.
Let us consider the path $\ga:[0,1]\to G$ from $1$ to
$$\bA\cdot\bB=\tau_{\infty,0}(\la_1)\cdot\tau_{\al,\be}(\la_2)$$
given by the suitably reparametrised product
$$\ga(t):=\tau_{\infty,0}(\la_1(t))\cdot\tau_{\al,\be}(\la_2(t))$$
with $\la_j(t)=1+t(\la_j-1)$ for $j=1,2$.
Let $\tilde\ga:[0,1]\to\Gm$ be the lift of this path covering with $\tilde\ga(0)=e$.
The path~$\ga$ is homotopic to the path
$$\de:=(\tau_{\infty,0}|_{[1,\la_1]})*(\tau_{\infty,0}(\la_1)\cdot\tau_{\al,\be}|_{[1,\la_2]}),$$
where $*$ means to go along the first path and then along the second path.
%(The homotopy between~$\ga$ and~$\de$ is given by $H(S,T)=\tau_{\infty,0}(\la_1^S)\cdot\tau_{\al,\be}(\la_2^T)$.)
It is clear that
$$\tilde\de:=([t_{\infty,0}]|_{[1,\la_1]})*([t_{\infty,0}(\la_1)]\cdot[t_{\al,\be}]|_{[1,\la_2]})$$
is the lift of the path~$\de$ with the starting point~$e$.
The end point of the lifted path~$\tilde\de$ is $[t_{\infty,0}(\la_1)\cdot t_{\al,\be}(\la_2)]$.
Since the path~$\ga$ is homotopic to~$\de$, the lift~$\tilde\ga$ has the same end point as~$\tilde\de$, hence
$$
  \tilde\ga(1)
  =[t_{\infty,0}(\la_1)\cdot t_{\al,\be}(\la_2)]
  =A\cdot B.
$$
So we have to compute
$$\levm(A\cdot B)-\levm(A)-\levm(B)=\levm(\tilde\ga(1)).$$

\myskip
On the other hand we have
\begin{align*}
  \ga(t)
  &=\tau_{\infty,0}(\la_1(t))\cdot\tau_{\al,\be}(\la_2(t))\\
  &=\left[\frac{1}{(\al-\be)\cdot\sqrt{\la_1(t)\la_2(t)}}\cdot
          \begin{pmatrix}\la_1(t)(\la_2(t)\al-\be)&-\la_1(t)(\la_2(t)-1)\al\be\\ \la_2(t)-1&\al-\la_2(t)\be\end{pmatrix}
    \right].
\end{align*}
Let $\Phi(t)=\varphi(\tilde\ga(t))$.
We obtain
$$\tan\frac{\Phi(t)}{2}=-\frac{(\la_1(t)\al\be+1)(\la_2(t)-1)}{(\la_1(t)\la_2(t)+1)\al-(\la_1(t)+\la_2(t))\be}.$$
We observe that $\Phi(0)=0$.
The denominator of this fraction is
\begin{align*}
  f(t)
  &:=(\la_1(t)\la_2(t)+1)\al-(\la_1(t)+\la_2(t))\be\\
  &=(\la_1(t)+\la_2(t))(\al-\be)+(\la_1(t)-1)(\la_2(t)-1)\al.
\end{align*}
Since $\la_1(t)+\la_2(t)>0$ and $(\la_1(t)-1)(\la_2(t)-1)\ge0$,
if $\al$ and $\al-\be$ are both positive or both negative,
\ie in the cases $\al<\be<0$, $\al<0<\be$, $\be<0<\al$ and $0<\be<\al$,
we have $f(t)\ne0$ for $t\in[0,1]$.
This implies $\Phi(1)\in(-\pi,\pi)$ and hence
$$\levm(A\cdot B)-\levm(A)-\levm(B)=\levm(\tilde\ga(1))=0.$$

\myskip
In the cases $0<\al<\be$ resp.\ $\be<\al<0$ we have to look carefully at the argument~$\Phi(t)$.
The argument $\Phi(t)$ satisfies the equation
$$\tan\frac{\Phi(t)}{2}=\frac{-(\la_1(t)\al\be+1)(\la_2(t)-1)}{f(t)}$$
where
\begin{align*}
  f(t)
  &=(\la_1(t)\la_2(t)+1)\al-(\la_1(t)+\la_2(t))\be\\
  &=\al(\la_1-1)(\la_2-1)\cdot t^2+(\al-\be)((\la_1-1)+(\la_2-1))\cdot t+2(\al-\be)
\end{align*}
is a quadratic function.
Since $\al\cdot\be>0$, $\la_1(t)>0$ and $\la_2(t)\ge1$,
the sign of $\tan(\Phi(t)/2)$ is opposite to the sign of~$f(t)$.

\myskip
Let us assume that $0<\al<\be$.
The coefficient $\al(\la_1-1)(\la_2-1)$ by~$t^2$ in~$f(t)$ is positive, hence the function $f$ is concave.
We observe that
\begin{align*}
  f(0)&=2(\al-\be)<0,\\
  f(1)&=(\la_1\la_2+1)\al-(\la_1+\la_2)\be.
\end{align*}
There are two cases, $f(1)\le0$ and $f(1)>0$.

\myskip
Let us assume that $0<\al<\be$ and $f(1)\le0$.
Then $f(0)<0$, $f(1)\le0$ and $f$ concave implies $f(t)<0$ for $t\in[0,1)$,
hence $\tan(\Phi(t)/2)>0$ for $t\in(0,1)$.
This implies $\Phi(t)\in[0,\pi]$ and hence
$$\levm(A\cdot B)-\levm(A)-\levm(B)=\levm(\tilde\ga(1))=0.$$

\myskip
Let us assume that $0<\al<\be$ and $f(1)>0$, which is equivalent to
$$0<\frac{\la_1+\la_2}{1+\la_1\la_2}\cdot\be<\al<\be.$$
Then $f(0)<0$, $f(1)>0$ and $f$ quadratic implies that there is $t_0\in(0,1)$
such that $f(t_0)=0$, $f(t)<0$ for $t\in[0,t_0)$, and $f(t)>0$ for $t\in(t_0,1]$.
Hence $\tan(\Phi(t)/2)>0$ for $t\in(0,t_0)$ and $\tan(\Phi(t)/2)<0$ for $t\in(t_0,1]$.
This implies $\Phi(1)\in(\pi,2\pi]$ and hence
$$\levm(A\cdot B)-\levm(A)-\levm(B)=\levm(\tilde\ga(1))=1.$$

\myskip
Let us assume that $\be<\al<0$.
The coefficient $\al(\la_1-1)(\la_2-1)$ by~$t^2$ in~$f(t)$ is negative, hence the function $f$ is convex.
We observe
\begin{align*}
  f(0)&=2(\al-\be)>0,\\
  f(1)&=(\la_1\la_2+1)\al-(\la_1+\la_2)\be.
\end{align*}
There are two cases, $f(1)>0$ and $f(1)\le0$.

\myskip
Let us assume that $\be<\al<0$ and $f(1)>0$.
Then $f(0)>0$, $f(1)>0$ and $f$ convex implies $f(t)\ne0$ for $t\in[0,1]$,
hence $\levm(A\cdot B)-\levm(A)-\levm(B)=0$ as before.

\myskip
Let us assume that $\be<\al<0$ and $f(1)\le0$, which is equivalent to
$$\be<\al\le\frac{\la_1+\la_2}{1+\la_1\la_2}\cdot\be<0.$$
Then $f(0)>0$, $f(1)\le0$ and $f$ quadratic implies that there is $t_0\in(0,1]$
such that $f(t_0)=0$, $f(t)>0$ for $t\in[0,t_0)$, and $f(t)<0$ for $t\in(t_0,1)$.
Hence $\tan(\Phi(t)/2)<0$ for $t\in(0,t_0)$ and $\tan(\Phi(t)/2)>0$ for $t\in(t_0,1)$.
This implies $\Phi(1)\in[-2\pi,-\pi]$ and hence
$$\levm(A\cdot B)-\levm(A)-\levm(B)=\levm(\tilde\ga(1))=-1.\qedhere$$
\end{proof}

The results similar to Lemma~\ref{lem-product-of-hyps} hold also for products
of hyperbolic and parabolic elements as well as for products of parabolic elements:

\begin{lem}
\label{lem-product-of-hyp-and-par}
For preimages $A$ resp.~$B$ in~$\Gm$
of the elements $\bA=\tau_{\infty,0}(\la_1)$ resp.\ $\bB=\pi_{\al}(\la_2)$
for some $\la_1>1$, $\la_2>0$ and $\al>0$
the difference
$$\levm(A\cdot B)-\levm(A)-\levm(B)$$
is equal

\begin{enumerate}[$\bullet$]
\item to~$0$ if
$$\la_2\al\le\frac{\la_1+1}{\la_1-1},$$
\item and to~$1$ otherwise.
\end{enumerate}
\end{lem}

\begin{lem}
\label{lem-product-of-pars}
For preimages $A$ resp.~$B$ in~$\Gm$
of the elements $\bA=\pi_{\infty}(\la_1)$ resp.\ $\bB=\pi_{\al}(\la_2)$
for some $\la_1,\la_2,\al>0$
the difference
$$\levm(A\cdot B)-\levm(A)-\levm(B)$$
is equal

\begin{enumerate}[$\bullet$]
\item to~$0$ if $\la_1\la_2\le2$,
\item and to~$1$ otherwise.
\end{enumerate}
\end{lem}

We omit the proofs of Lemmata~\ref{lem-product-of-hyp-and-par} and~\ref{lem-product-of-pars},
which are along the lines of the proof of Lemma~\ref{lem-product-of-hyps}.

\section{Higher spin structures and lifts of Fuchsian groups}

\label{sec-spins-and-lifts}

\subsection{Higher spin structures}

Let $E\to P$ be complex line bundle over a hyperbolic Riemann surface $P$.
Let $\lat$ be a torsionfree Fuchsian group such that $P=\hyp/\lat$.
Let $L\to\hyp$ be the induced complex line bundle over $\hyp$.
Let $L\simeq\hyp\times\c$ be a trivialization of the bundle $L$.
With respect to this trivialization the action of $\lat$ on $L$ is given by
$$g\cdot(z,t)=(g(z),\de(g,z)\cdot t),$$
where $\de:\lat\times\h\to\c^*$ is a map such that the function
$\de_g=\de|_{\{g\}\times\hyp}$ is holomorphic for any $g\in\lat$ and for
any $g_1,g_2\in\lat$ we have
$$\de_{g_2\cdot g_1}=(\de_{g_2}\circ g_1)\cdot\de_{g_1}.$$
The map $\de$ is called the {\it transition map} of the bundle $E\to P$ with
respect to the given trivialization.

\myskip
In particular, if $E$ is the cotangent bundle of the surface $P$, then the transition map can be chosen so that $\de_g=(g')^{-1}$.
If $E$ is the tangent bundle of the surface $P$, then the transition map can be chosen so that $\de_g=g'$.
Let $E_1\to P$, $E_2\to P$ be two complex line bundles over a Riemann surface $P$,
and let $\de_1$ resp.\ $\de_2$ be their transition maps,
then $\de_1\cdot\de_2$ is a transition map of the bundle $E_1\otimes E_2\to P$.
In particular, if $\de$ is the transition map of the bundle $E\to P$,
then $\de^m$ is a transition map of the bundle $E^m=E\otimes\cdots\otimes E\to P$
(with respect to the induced trivialization).

\myskip
An {\it $m$-spin structure} on a Riemann surface~$P$
is a transition map $\de$ of a complex line bundle $E\to P$
that satisfies the condition $\de_g^m=(g')^{-1}$,
\ie the induced transition map $\de^m$ of the bundle $E^m\to P$ coincides with the transition map of the cotangent bundle of~$P$.

\begin{rem}
An {\it $m$-co-spin structure} on a Riemann surface~$P$
is a transition map $\de$ of a complex line bundle $E\to P$
that satisfies the condition $\de_g^m=g'$,
\ie the induced transition map $\de^m$ of the bundle $E^m\to P$ coincides with the transition map of the tangent bundle of~$P$.
There is a one-to-one correspondence between $m$-spin and $m$-co-spin
structures on a Riemann surface given by taking $E$ to~$E^{-1}$ resp., in terms of transition maps, by taking $\de$ to $\de^{-1}$.
In the following we consider the $m$-spin structures.
\end{rem}

\begin{rem}
A complex line bundle $E\to P$ is said to be {\it $m$-spin} if the bundle $E^m\to P$ is isomorphic to the cotangent bundle of~$P$.
For a compact Riemann surface $P$ there is a 1-1-correspondence between $m$-spin structures on~$P$ and $m$-spin bundles over~$P$.
\end{rem}

\begin{mydef}
A {\it lift} of the Fuchsian group $\lat$ into $\Gm$
is a subgroup $\lats$ of $\Gm$ such that the restriction
of the covering map $\Gm\to G$ to $\lats$ is an isomorphism
between $\lats$ and $\lat$.
\end{mydef}

\begin{prop}
\label{spinstructures-lifts}
There is a 1-1-correspondence between $m$-spin structures on the
Riemann surface $P=\hyp/\lat$ and lifts of $\lat$ into $\Gm$.
\end{prop}

\begin{proof}
We use the description in Proposition~\ref{fract-autforms} of the covering $\Gm$ as the group
$$\{(g,\de)\in G\times\Hol(\hyp,\c^*)\st \de^m(z)=g'(z)~{\rm for~all}~z\in\hyp\}$$
with multiplication
$$(g_2,\de_2)\cdot(g_1,\de_1)=(g_2\cdot g_1,(\de_2\circ g_1)\cdot\de_1).$$
From this description of~$\Gm$ we see that chosing a lift of~$\lat$ into~$\Gm$
amounts to chosing for each~$g\in\lat$ a holomorphic function~$\de_g:\hyp\to\c^*$ that satisfies the condition~$\de_g^m=g'$
so that that the family~$\{\de_g\}_{g\in\lat}$ satisfies the condition
$$\de_{g_2\cdot g_1}=(\de_{g_2}\circ g_1)\cdot\de_{g_1}\quad\text{for any}~g_1,g_2\in\lat.$$

%On the one hand using the description of the covering $\Gm$ from Proposition~\ref{fract-autforms}
% we see that there is a 1-1-correspondence between the lifts of $\lat$ into $\Gm$
% and the families $\{\de_g\}_{g\in\lat}$ of holomorphic functions $\hyp\to\c^*$ such that
%for any $g\in\lat$
%$$\de_g^m=g'$$
%and for any $g_1,g_2\in\lat$
%$$\de_{g_2\cdot g_1}=(\de_{g_2}\circ g_1)\cdot\de_{g_1}.$$

On the other hand as we explained in this section an $m$-spin structure on the Riemann surface~$P=\hyp/\lat$ is described
by a family~$\{\de_g\}_{g\in\lat}$ that satisfies the same conditions
$$\de_g^m=g'\quad\text{for any}~g\in\lat$$
and
$$\de_{g_2\cdot g_1}=(\de_{g_2}\circ g_1)\cdot\de_{g_1}\quad\text{for any}~g_1,g_2\in\lat.$$
Hence there is a 1-1-correspondence between $m$-spin structures on the
Riemann surface $P=\hyp/\lat$ and lifts of $\lat$ into $\Gm$.
\end{proof}

\subsection{Finitely generated Fuchsian groups}

\myskip
In this section we are going to describe finitely generated Fuchsian groups using standard sets of generators.
Definition~\ref{def-short-seq-set} corresponds to the case of the fundamental group of a sphere
with punctures and holes with the total number of punctures and holes equal to three.
Definition~\ref{def-long-seq-set} corresponds to the case of the fundamental group of a sphere with an arbitrary number of punctures and holes.
Finally in Definition~\ref{def-genus-seq-set} we describe the fundamental group of a surface (of possibly higher genus~$g$) with punctures and holes.
Hereby we use the observation that if we think of the surface as a sphere with $g$~handles (and punctures and holes) and cut all handles along their waist curves
then we obtain a sphere with the same number of punctures and the number of holes increased by~$2g$.

\begin{mydef}
For two elements~$C_1$ and~$C_2$ in~$G$ with finite fixed points in~$\r$ we say that $C_1<C_2$
if all fixed point of~$C_1$ are smaller than any fixed point of~$C_2$.
\end{mydef}

\begin{figure}
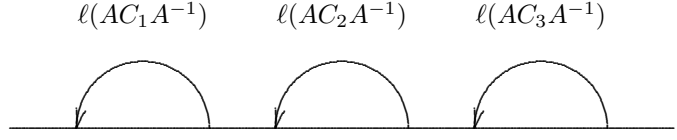

  \begin{center}
    \forcehmode
      \bgroup
        \beginpicture
          \setcoordinatesystem units <25 bp,25 bp>
          \multiput {\phantom{$\bullet$}} at -5 -1 5 2 /
          \circulararc 180 degrees from -2 0 center at -3 0
          \arrow <7pt> [0.2,0.5] from -3.99 0.04 to -4 0
          \put {$\ell(AC_1A^{-1})$} [b] <0pt,\baselineskip> at -3 1
          \circulararc 180 degrees from  1 0 center at  0 0
          \arrow <7pt> [0.2,0.5] from -0.99 0.04 to -1 0
          \put {$\ell(AC_2A^{-1})$} [b] <0pt,\baselineskip> at 0 1
          \circulararc 180 degrees from  4 0 center at  3 0
          \arrow <7pt> [0.2,0.5] from  2.01 0.04 to  2 0
          \put {$\ell(AC_3A^{-1})$} [b] <0pt,\baselineskip> at 3 1
          \plot -5 0 5 0 /
        \endpicture
      \egroup
  \end{center}
  \caption{Axes of a sequential set of type $(0,3,0)$}
  \label{fig-axes-seqset}
\end{figure}

\begin{mydef}
\label{def-short-seq-set}
A {\it sequential set of type}~$(0,l_h,l_p)$ with $l_h+l_p=3$ is a triple of elements $(C_1,C_2,C_3)$ in~$G$
such that
the elements $C_i$ for~$1\le i\le l_h$ are hyperbolic,
the elements $C_i$ for~$l_h<i\le3$ are parabolic,
their product is
$$C_1\cdot C_2\cdot C_3=1,$$
and for some element $A\in G$ the elements~$\{\tilde C_i=A C_i A^{-1}\}_{i=1,2,3}$ are positive, have finite fixed points
and satisfy~$\tilde C_1<\tilde C_2<\tilde C_3$.
(Figure~\ref{fig-axes-seqset} illustrates the position of the axes of the elements~$\tilde C_i$ for a sequential set of type $(0,3,0)$, \ie when all elements are hyperbolic.)
%It is clear what the similar picture looks like in presence of parabolic elements.)
\end{mydef}

% Picture

\begin{mydef}
\label{def-long-seq-set}
A {\it sequential set of type}~$(0,l_h,l_p)$ is an $n$-tuple of elements
$$(C_1,\dots,C_n)$$
with $n=l_h+l_p$ in~$G$ such that
the elements $C_1,\dots,C_{l_h}$ are hyperbolic,
the elements $C_{l_h+1},\dots,C_n$ are parabolic,
and for any $j\in\{2,\dots,n-1\}$ the triple
$$(C_1\cdots C_{j-1},C_j,C_{j+1}\cdots C_n)$$
is a sequential set (of type $(0,3,0)$, $(0,2,1)$, $(0,1,2)$ or~$(0,0,3)$).
\end{mydef}

\begin{mydef}
\label{def-genus-seq-set}
A {\it sequential set of type}~$(g,l_h,l_p)$ is a $(2g+l_h+l_p)$-tuple of elements
$$(A_1,\dots,A_g,B_1,\dots,B_g,C_{g+1},\dots,C_{g+l_h+l_p})$$
%with $n=l_h+l_p$
in~$G$ such that
the elements $A_1,\dots,A_g,B_1,\dots,B_g$ and~$C_{g+1},\dots,C_{g+l_h}$ are hyperbolic,
the elements $C_{g+l_h+1},\dots,C_{g+l_h+l_p}$ are parabolic,
and the tuple
$$(A_1,B_1A_1^{-1}B_1^{-1},\dots,A_g,B_gA_g^{-1}B_g^{-1},C_{g+1},\dots,C_{g+l_h+l_p})$$
is a sequential set of type $(0,2g+l_h,l_p)$.
\end{mydef}

\begin{mydef}
We call a Riemann surface of genus~$g$ with $l_h$ holes and $l_p$ punctures
a {\it Riemann surface of type} $(g,l_h,l_p)$.
\end{mydef}

\begin{figure}
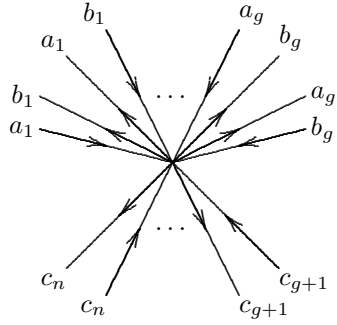

  \begin{center}
    \forcehmode
      \bgroup
        \beginpicture
          \setcoordinatesystem units <25 bp,25 bp>
          \multiput {\phantom{$\bullet$}} at -2 -2 2 2 /
          \plot 0 0 -2 0.5 /
          \arrow <7pt> [0.2,0.5] from -2 0.5 to -1 0.25
          \plot 0 0 -2 1 /
          \arrow <7pt> [0.2,0.5] from 0 0 to -1 0.5
          \plot 0 0 -1.6 1.6 /
          \arrow <7pt> [0.2,0.5] from 0 0 to -0.8 0.8
          \plot 0 0 -1 2 /
          \arrow <7pt> [0.2,0.5] from -1 2 to -0.5 1
          \plot 0 0 1 2 /
          \arrow <7pt> [0.2,0.5] from 1 2 to 0.5 1
          \plot 0 0 1.6 1.6 /
          \arrow <7pt> [0.2,0.5] from 0 0 to 0.8 0.8
          \plot 0 0 2 1 /
          \arrow <7pt> [0.2,0.5] from 0 0 to 1 0.5
          \plot 0 0 2 0.5 /
          \arrow <7pt> [0.2,0.5] from 2 0.5 to 1 0.25
          \plot 0 0 1.6 -1.6 /
          \arrow <7pt> [0.2,0.5] from 1.6 -1.6 to 0.8 -0.8
          \plot 0 0 1 -2 /
          \arrow <7pt> [0.2,0.5] from 0 0 to 0.5 -1
          \plot 0 0 -1 -2 /
          \arrow <7pt> [0.2,0.5] from -1 -2 to -0.5 -1
          \plot 0 0 -1.6 -1.6 /
          \arrow <7pt> [0.2,0.5] from 0 0 to -0.8 -0.8
          \put {$a_1$} [r] <-2pt,0pt> at -2 0.5
          \put {$b_1$} [r] <-2pt,0pt> at -2 1
          \put {$a_1$} [br] <0pt,2pt> at -1.6 1.6
          \put {$b_1$} [br] <0pt,2pt> at -1 2
          \put {$a_g$} [bl] <0pt,2pt> at 1 2
          \put {$b_g$} [bl] <0pt,2pt> at 1.6 1.6
          \put {$a_g$} [l] <2pt,0pt> at 2 1
          \put {$b_g$} [l] <2pt,0pt> at 2 0.5
          \put {$c_{g+1}$} [tl] <0pt,-2pt> at 1 -2
          \put {$c_{g+1}$} [tl] <0pt,-2pt> at 1.6 -1.6
          \put {$c_n$} [tr] <0pt,-2pt> at -1.6 -1.6
          \put {$c_n$} [tr] <0pt,-2pt> at -1 -2
          \put {$\dots$} at 0 1
          \put {$\dots$} at 0 -1
        \endpicture
      \egroup
  \end{center}
  \caption{Standard basis}
  \label{fig-basis}
\end{figure}

\begin{mydef}
\label{def-standard-basis}
We define the product $ab$ of two contours $a$ and~$b$ in $\pi_1(P,p)$
as the contour given by the path of~$b$ followed by the path of~$a$.
A {\it standard basis} of a fundamental group $\pi_1(P,p)$
of a surface $P$ of type $(g,l_h,l_p)$ is a set of generators
$$\{a_i,b_i~(i=1,\dots,g),c_i~(i=g+1,\dots,n)\}$$
with a single defining relation
$$\prod_{i=1}^g [a_i,b_i] \prod_{i=g+1}^n c_i=1$$
and represented by a set of simple contours
$$\{\tilde a_i,\tilde b_i(i=1,\dots,g),\tilde c_i(i=g+1,\dots,n)\}$$
with the following properties:
\begin{enumerate}[1)]
\item
the contour $\tilde c_i$ encloses a hole in $P$ for $i=g+1,\dots,g+l_h$ and a puncture for $i=g+l_h+1,\dots,n$,
\item
$\tilde a_i\cap\tilde b_j=\tilde a_i\cap\tilde c_j=\tilde b_i\cap\tilde c_j=\tilde c_i\cap\tilde c_j=\{p\}$.
\item
in a neighbourhood of the point $p$, the contours are placed as is shown in
Figure~\ref{fig-basis}.
\end{enumerate}
\end{mydef}

% Picture

The relation between sequential sets and Fuchsian groups is exploited
in \cite{Nbook} (Chapter~1, Theorem~1.1, Lemma~2.1, Theorem~2.1).
We recall here the results:

\begin{thm}
\label{seqset-fuchgr}
A sequential set $V$ of type $(g,l_h,l_p)$ generates a Fuchsian group~$\lat$
such that the surface $P=\hyp/\lat$ is of type $(g,l_h,l_p)$.
The isomorphism $\Phi:\lat\to\pi_1(P,p)$, induced by the natural projection $\Psi:\hyp\to P$,
maps the sequential set $V$ to a standard basis of $\pi_1(P,p)$.
\end{thm}

\begin{thm}
\label{thm2.1inN}
Let $\lat$ be a Fuchsian group such that the surface $P=\hyp/\lat$ is of type $(g,l_h,l_p)$.
Let $p$ be a point in~$P$.
Let $\Psi:\hyp\to P$ be the natural projection.
Choose $q\in\Psi^{-1}(p)$ and let $\Phi:\lat\to\pi_1(P,p)$ be the induced isomorphism.
Let
$$v=\{a_i,b_i~(i=1,\dots,g),c_i~(i=g+1,\dots,n)\}$$
be a standard basis of $\pi_1(P,p)$.
In this case,
\begin{align*}
  V=\Phi^{-1}(v)&=\{\Phi^{-1}(a_i),\Phi^{-1}(b_i)~(i=1,\dots,g),\Phi^{-1}(c_i)~(i=g+1,\dots,n)\}\\
                &=\{A_i,B_i~(i=1,\dots,g),C_i~(i=g+1,\dots,n)\}
\end{align*}
is a sequential set of type $(g,l_h,l_p)$.
\end{thm}

Now we recall the classification of free Fuchsian groups of rank~$2$ (see \cite{Nbook}, Chapter~1, Lemma~3.2, 3.3, 3.4):

\begin{lem}
\label{lem-FG1}
The set
$(C_1=\tau_{\infty,0}(\la_1),C_2=\tau_{\al,\be}(\la_2),C_3)$
with $\la_1,\la_2>1$ is a sequential set of type $(0,3,0)$ or $(0,2,1)$ if and only if
$$C_3=(C_1\cdot C_2)^{-1}$$
and
$$
  0<\left(\frac{\sqrt{\la_1}+\sqrt{\la_2}}{1+\sqrt{\la_1\la_2}}\right)^{\text{\rlap{2}}}\be
  \le\al<\be<\infty.
$$
Then the set $(C_1,C_2,C_3)$ is of type $(0,2,1)$,
\ie the element~$C_3$ is parabolic, if and only if
$$\al=\left(\frac{\sqrt{\la_1}+\sqrt{\la_2}}{1+\sqrt{\la_1\la_2}}\right)^2\be.$$
\end{lem}

\begin{lem}
\label{lem-FG2}
The set
$(C_1=\tau_{\infty,0}(\la_1),C_2=\pi_{\al}(\la_2),C_3)$
with $\la_1>1$, $\la_2>0$ is a sequential set of type $(0,2,1)$ or $(0,1,2)$ if and only if
$$C_3=(C_1\cdot C_2)^{-1}$$
and
$$\la_2\al\ge\frac{\sqrt{\la_1}+1}{\sqrt{\la_1}-1}.$$
Then the set $(C_1,C_2,C_3)$ is of type $(0,1,2)$,
\ie the element~$C_3$ is parabolic,
if and only if the last inequality is an equality.
\end{lem}

\begin{lem}
\label{lem-FG3}
If the set~$(A_1,B_1,C_1)$ is a sequential set of type $(1,1,0)$
then the axes of the hyperbolic elements $A_1$ and~$B_1$ intersect.
\end{lem}

\subsection{Lifting sets of generators of Fuchsian groups}

\label{first-rem-FG-into-tG}

\myskip
In this subsection let us denote by $[\cdot]$ the image of an element in~$\tG$ under the covering map $\tG\to\Gm$.

\begin{lem}
\label{lem-liftonerel}
Let $\lat$ be a Fuchsian group of type $(g,l_h,l_p)$ generated by the sequential set
\begin{align*}
  \bV
  &=\{\bA_i,\bB_i~(i=1,\dots,g),\bC_i~(i=g+1,\dots,n)\}\\
  &=\{\bD_j~(j=1,\dots,n+g)\}.
\end{align*}
Let
\begin{align*}
  V&=\{A_i,B_i~(i=1,\dots,g),C_i~(i=g+1,\dots,n)\}\\
  &=\{D_j~(j=1,\dots,n+g)\}
\end{align*}
be a set of the lifts of the elements of the sequential set $\bV$ into~$\Gm$,
\ie the image of $D_j$ in $G$ is $\bD_j$.
Then the subgroup $\lats$ of $\Gm$ generated by $V$ is a lift of $\lat$ into~$\Gm$ if and only if
$$\prod\limits_{i=1}^g\,[A_i,B_i]\cdot\prod\limits_{i=g+1}^n\,C_i=e.$$
\end{lem}

\begin{proof}
For any choice of the set of lifts $V$ the restriction of the covering map
$\Gm\to G$ to the group $\lats$ generated by $V$ is a homomorphism with image $\lat$.
There is only one relation
$$\prod\limits_{i=1}^g\,[\bA_i,\bB_i]\cdot\prod\limits_{i=g+1}^n\,\bC_i=e$$
in $\lat$,
hence the equality
$$\prod\limits_{i=1}^g\,[A_i,B_i]\cdot\prod\limits_{i=g+1}^n\,C_i=e$$
ensures injectivity of this homomorphism.
\end{proof}

%\subsection{Lifts of free Fuchsian groups of rank $2$}
%\label{FFG2-lifts}

\begin{lem}
\label{lem-lift-030}
Let $(C_1,C_2,C_3)$ be a triple of elements in~$\Gm$
such that their images $(\bC_1,\bC_2,\bC_3)$ in~$G$ form a sequential set of type~$(0,l_h,l_p)$ with $l_h+l_p=3$.
Then $C_1\cdot C_2\cdot C_3=e$ if and only if
$$\levm(C_1)+\levm(C_2)+\levm(C_3)=-1.$$
Moreover
$$\levm(C_1\cdot C_2)=\levm(C_1)+\levm(C_2)+1.$$
\end{lem}

\begin{proof}
We first prove $\levm(C_1\cdot C_2)=\levm(C_1)+\levm(C_2)+1$ separately
for sequential sets of types $(0,3,0)$ and $(0,2,1)$ using Lemma~\ref{lem-product-of-hyps},
for sequential sets of types $(0,1,2)$ and $(0,2,1)$ using Lemma~\ref{lem-product-of-hyp-and-par},
and for sequential sets of type~$(0,0,3)$ using Lemma~\ref{lem-product-of-pars}.
\begin{enumerate}[$\bullet$]
\item
We first assume that the elements $C_1$ and $C_2$ are hyperbolic.
Up to conjugation we can assume that $\bC_1=\tau_{\infty,0}(\la_1)$ and $\bC_2=\tau_{\al,\be}(\la_2)$
for some $\la_1,\la_2>1$ and $\al,\be\in\r\bs\{0\}$.
Lemma~\ref{lem-FG1} implies that $\al,\be>0$ and
$$
  1>\frac{\al}{\be}
  >\left(\frac{\sqrt{\la_1}+\sqrt{\la_2}}{1+\sqrt{\la_1\la_2}}\right)^2
  =\frac{\la_1+\la_2+2\sqrt{\la_1\la_2}}{1+\la_1\la_2+2\sqrt{\la_1\la_2}}
  >\frac{\la_1+\la_2}{1+\la_1\la_2}.
$$
According to Lemma~\ref{lem-product-of-hyps} the inequalities $\al\cdot\be>0$ and
$$\frac{\la_1+\la_2}{1+\la_1\la_2}\cdot\be<\al<\be$$
imply
$$\levm(C_1\cdot C_2)=\levm(C_1)+\levm(C_2)+1.$$
\item
We now assume that the element~$C_1$ is hyperbolic and the element~$C_2$ is parabolic.
Up to conjugation we can assume that $\bC_1=\tau_{\infty,0}(\la_1)$ and $\bC_2=\pi_{\al}(\la_2)$
for some $\la_1>1$, $\la_2>0$, and $\al>0$.
Lemma~\ref{lem-FG2} implies that
$$
  \la_2\al
  \ge\frac{\sqrt{\la_1}+1}{\sqrt{\la_1}-1}
  =\frac{(\sqrt{\la_1}+1)^2}{(\sqrt{\la_1}+1)\cdot(\sqrt{\la_1}-1)}
  =\frac{\la_1+1+2\sqrt{\la_1}}{\la_1-1}
  >\frac{\la_1+1}{\la_1-1}.
$$
According to Lemma~\ref{lem-product-of-hyp-and-par} the inequalities $\al>0$ and
$$\la_2\al>\frac{\la_1+1}{\la_1-1}$$
imply
$$\levm(C_1\cdot C_2)=\levm(C_1)+\levm(C_2)+1.$$
\item
We now assume that the elements~$C_1$ and~$C_2$ are parabolic.
Up to conjugation we can assume that
$$
  \bC_1=\pi_{\infty}(1)=\left[\begin{pmatrix} 1&1\\ 0&1\end{pmatrix}\right]
  \quad\hbox{and}\quad
  \bC_2=\pi_1(4)=\left[\begin{pmatrix} -3&4\\ -4&5\end{pmatrix}\right].
$$
Then according to Lemma~\ref{lem-product-of-pars}
$$\levm(C_1\cdot C_2)=\levm(C_1)+\levm(C_2)+1.$$
\item
In all three cases we proved that $\levm(C_1\cdot C_2)=\levm(C_1)+\levm(C_2)+1$.
For the inverse element this implies
$$\levm((C_1\cdot C_2)^{-1})=-\levm(C_1\cdot C_2)=-(\levm(C_1)+\levm(C_2)+1).$$
The image of the element $C_1\cdot C_2\cdot C_3$ is $\bC_1\cdot\bC_2\cdot\bC_3=1$,
hence $C_1\cdot C_2\cdot C_3=[u]^l$ for some $l\in\z$.
This implies
$$\levm(C_3)=\levm((C_1\cdot C_2)^{-1}\cdot u^l)=l-(\levm(C_1)+\levm(C_2)+1),$$
hence we obtain that $C_1\cdot C_2\cdot C_3=e$ if and only if $l=0$,
\ie if and only if
$$\levm(C_1)+\levm(C_2)+\levm(C_3)=-1.\qedhere$$
\end{enumerate}
\end{proof}

\begin{lem}
\label{lem-lift-110}
Let $(A_1,B_1,C_1)$ be a triple of elements in~$\Gm$
such that there images $(\bA_1,\bB_1,\bC_1)$ in~$G$ form a sequential set of type~$(1,1,0)$.
Then
$$[A_1,B_1]\cdot C_1=e$$
if and only if
$$\levm(C_1)=-1.$$
Moreover,
$$\levm(A_1\cdot B_1)=\levm(A_1)+\levm(B_1).$$
\end{lem}

\begin{proof}
Lemma~\ref{lem-FG3} implies that the axes of the hyperbolic elements~$A_1$ and~$B_1$ intersect.
According to Corollary~\ref{cor-product-of-crossing-hyps} this implies
$$\levm(A_1\cdot B_1)=\levm(A_1)+\levm(B_1).$$
The triple $(\bA_1,\bB_1,\bC_1)$ is a sequential set of type~$(1,1,0)$.
By the definition of sequential sets the triple
$(\bA_1,\bB_1\bA_1^{-1}\bB_1^{-1},\bC_1)$ is a sequential set of type~$(0,3,0)$.
According to Lemma~\ref{lem-lift-030} we have
$$[A_1,B_1]\cdot C_1=A_1\cdot(B_1A_1^{-1}B_1^{-1})\cdot C_1=e$$
if and only if
$$\levm(A_1)+\levm(B_1A_1^{-1}B_1^{-1})+\levm(C_1)=-1.$$
Since by Lemma~\ref{lem-conj}
$$\levm(B_1A_1^{-1}B_1^{-1})=\levm(A_1^{-1})=-\levm(A_1),$$
the last condition is equivalent to
$$\levm(C_1)=-1.\qedhere$$
\end{proof}

\begin{lem}
\label{lem-lift-glm}
Let
$$V=\{A_i,B_i~(i=1,\dots,g),C_i~(i=g+1,\dots,n)\}$$
be a tuple of elements in~$\Gm$
such that the image
$$\bV=\{\bA_i,\bB_i~(i=1,\dots,g),\bC_i~(i=g+1,\dots,n)\}$$
in $G$ form a sequential set.
Then we have
$$\prod\limits_{i=1}^g\,[A_i,B_i]\cdot\prod\limits_{i=g+1}^n\,C_i=e$$
if and only if
$$\sum\limits_{i=g+1}^n\,\levm(C_i)=(2-2g)-(n-g).$$
(In the case $n=g$ this means $2-2g=0~\mod~m$.)
\end{lem}

\begin{proof}
We discuss the case $g=0$ first, and then we reduce the general case to the case $g=0$.
\begin{enumerate}[$\bullet$]
\item
Let $g=0$.
We prove that the statement is true for lifts of sequential sets of type $(0,l_h,l_p)$ by induction on~$l_h+l_p$.
The case $l_h+l_p=3$ is covered by Lemma~\ref{lem-lift-030}.
% Property~\ref{arf-prop-seqset} of $m$-Arf functions.
Assume that the statement is true for $l_h+l_p\le n-1$ and consider the case $l_h+l_p=n$.
By the definition of sequential sets the set $(\bC_1\cdot\bC_2,\bC_3,\dots,\bC_n)$ is a sequential set.
Hence by our assumption $(C_1\cdot C_2)\cdot C_3\cdots C_n=e$ if and only if
$$\levm(C_1\cdot C_2)+\levm(C_3)+\cdots+\levm(C_n)=2-(n-1)=(2-n)+1.$$
Moreover, by the definition of sequential sets the set $(\bC_1,\bC_2,\bC_3\cdots\bC_n)$ is a sequential set too,
hence by Lemma~\ref{lem-lift-030}
% Property~\ref{arf-prop-seqset} of $m$-Arf functions
we have
$$\levm(C_1\cdot C_2)=\levm(C_1)+\levm(C_2)+1.$$
The last two equations imply that
$$C_1\cdots C_n=e$$
if and only if
$$\levm(C_1)+\cdots+\levm(C_n)=2-n.$$
\item
We now consider the general case.
By the definition of sequential sets the set
$$(\bA_1,\bB_1\bA_1^{-1}\bB_1^{-1},\dots,\bA_g,\bB_g\bA_g^{-1}\bB_g^{-1},\bC_1,\dots,\bC_n)$$
is a sequential set of type $(0,2g+l_h,l_p)$, hence
$$
  \prod\limits_{i=1}^g\,[A_i,B_i]\cdot\prod\limits_{i=g+1}^n C_i
  =\prod\limits_{i=1}^g\,(A_i\cdot B_i A_i^{-1}B_i^{-1})\cdot\prod\limits_{i=g+1}^n C_i
  =e
$$
if and only if
\begin{align*}
  \sum\limits_{i=1}^g\,(\levm(A_i)+\levm(B_i A_i^{-1}B_i^{-1}))+\sum\limits_{i=g+1}^n\levm(C_i)
  &=2-(n+g)\\
  &=(2-2g)-(n-g).
\end{align*}
From Lemma~\ref{lem-conj}
% Properties~\ref{arf-prop-conj} and~\ref{arf-prop-inv} of $m$-Arf functions
we obtain that
$$
  \levm(B_i A_i^{-1}B_i^{-1})
  =\levm(A_i^{-1})
  =-\levm(A_i),
$$
and hence
$$\levm(A_i)+\levm(B_i A_i^{-1}B_i^{-1})=0.\qedhere$$
\end{enumerate}
\end{proof}

\section{Higher Arf functions}

\label{sec-m-arf}

Let $\lat$ be a Fuchsian group
and $P=\hyp/\lat$ of type $(g,l_h,l_p)$ with $l_h+l_p=n-g$.
Let $p\in P$.
Let $\Psi:\hyp\to P$ be the natural projection.
Choose $q\in\Psi^{-1}(p)$ and let $\Phi:\lat\to\pi_1(P,p)$ be the induced isomorphism.

\subsection{Definition of higher Arf functions}

\label{m-arf}

\myskip
Let $\lats$ be a lift of $\lat$ in $\Gm$.

\begin{mydef}
\label{def-hsi-to-lats}
Let us consider a function $\hsi_{\lats}:\pi_1(P,p)\to\zm$ such that the following diagram commutes
$$
  \begin{CD}
   \lat            @>{\cong}>> \lats         \\
   @V{\Phi}VV         @VV{\levm|_{\lats}}V \\
   \pi_1(P,p)     @>{\hsi_{\lats}}>> \zm          \\
  \end{CD}
$$
% As for the function~$\levm$, all equations involving~$\hsi_{\lats}$ are to be understood as equations in~$\zm$.
\end{mydef}

\begin{lem}
\label{lem-hsi-rules}
Let $\al$, $\be$, and $\ga$ be simple contours in $P$
intersecting pairwise in exactly one point~$p$.
Let $a$, $b$, and $c$ be the corresponding elements of $\pi_1(P,p)$.
We assume that $a$, $b$, and $c$ satisfy the relations $a,b,c\ne1$ and $abc=1$.
Let $\<\cdot,\cdot\>$ be the intersection form on $\pi_1(P,p)$.
Then for $\hsi:=\hsi_{\lats}$
\begin{enumerate}[1.]
\item
$\hsi(ab)=\hsi(a)+\hsi(b)$
if the elements~$a$ and~$b$ can be represented by a pair of simple contours in $P$
intersecting in exactly one point~$p$
with $\<a,b\>\ne0$,
\item
$\hsi(ab)=\hsi(a)+\hsi(b)+1$
if the elements~$a$ and~$b$ can be represented by a pair of simple contours in $P$
intersecting in exactly one point~$p$
with $\<a,b\>=0$
and placed in a neighbourhood of the point~$p$ as shown in Figure~\ref{fig-pos-pair}.

% Picture

\begin{figure}
  \begin{center}
    \forcehmode
      \bgroup
        \beginpicture
          \setcoordinatesystem units <25 bp,25 bp>
          \multiput {\phantom{$\bullet$}} at -2 0 2 2 /
          \plot 0 0 -2 1 /
          \arrow <7pt> [0.2,0.5] from -2 1 to -1 0.5
          \plot 0 0 -1 2 /
          \arrow <7pt> [0.2,0.5] from 0 0 to -0.5 1
          \plot 0 0 1 2 /
          \arrow <7pt> [0.2,0.5] from 1 2 to 0.5 1
          \plot 0 0 2 1 /
          \arrow <7pt> [0.2,0.5] from 0 0 to 1 0.5
          \put {$a$} [br] <0pt,2pt> at -2 1
          \put {$a$} [br] <0pt,2pt> at -1 2
          \put {$b$} [bl] <0pt,2pt> at 2 1
          \put {$b$} [bl] <0pt,2pt> at 1 2
        \endpicture
      \egroup
  \end{center}
  \caption{$\hsi(ab)=\hsi(a)+\hsi(b)+1$}
  \label{fig-pos-pair}
\end{figure}

\item
$\hsi(ab)=\hsi(a)+\hsi(b)-1$
if the elements~$a$ and~$b$ can be represented by a pair of simple contours in $P$
intersecting in exactly one point~$p$
with $\<a,b\>=0$
and placed in a neighbourhood of the point~$p$ as shown in Figure~\ref{fig-neg-pair}.

% Picture

\begin{figure}
  \begin{center}
    \forcehmode
      \bgroup
        \beginpicture
          \setcoordinatesystem units <25 bp,25 bp>
          \multiput {\phantom{$\bullet$}} at -2 0 2 2 /
          \plot 0 0 -2 1 /
          \arrow <7pt> [0.2,0.5] from 0 0 to -1 0.5
          \plot 0 0 -1 2 /
          \arrow <7pt> [0.2,0.5] from -1 2 to -0.5 1
          \plot 0 0 1 2 /
          \arrow <7pt> [0.2,0.5] from 0 0 to 0.5 1
          \plot 0 0 2 1 /
          \arrow <7pt> [0.2,0.5] from 2 1 to 1 0.5
          \put {$a$} [br] <0pt,2pt> at -2 1
          \put {$a$} [br] <0pt,2pt> at -1 2
          \put {$b$} [bl] <0pt,2pt> at 2 1
          \put {$b$} [bl] <0pt,2pt> at 1 2
        \endpicture
      \egroup
  \end{center}
  \caption{$\hsi(ab)=\hsi(a)+\hsi(b)-1$}
  \label{fig-neg-pair}
\end{figure}

\item
for any standard basis
$$v=\{a_i,b_i~(i=1,\dots,g),c_i~(i=g+1,\dots,n)\}$$
of $\pi_1(P,p)$ we have
$$\sum\limits_{i=g+1}^n\,\hsi(c_i)=(2-2g)-(n-g).$$
\end{enumerate}
\end{lem}

\begin{proof}
According to Theorem~\ref{thm2.1inN}
either the set
$$V:=\{\Phi^{-1}(a),\Phi^{-1}(b),\Phi^{-1}(c)\}$$
or the set
$$V^{-1}:=\{\Phi^{-1}(a^{-1}),\Phi^{-1}(b^{-1}),\Phi^{-1}(c^{-1})\}$$
is sequential.
This sequential set can be of type $(0,k,m)$ (with $k+m=3$) or of type $(1,1,0)$.
If $V$ is a sequential set of type $(1,1,0)$, then according to Lemma~\ref{lem-lift-110} we obtain
$$\hsi(ab)=\hsi(a)+\hsi(b).$$
If $V$ is a sequential set of type $(0,k,m)$, then according to Lemma~\ref{lem-lift-030} we obtain
$$\hsi(ab)=\hsi(a)+\hsi(b)+1.$$
If $V^{-1}$ is a sequential set of type $(0,k,m)$, then according to Lemma~\ref{lem-lift-030} we obtain
$$\hsi(b^{-1}a^{-1})=\hsi(a^{-1})+\hsi(b^{-1})+1$$
and hence
\begin{align*}
  -\hsi(ab)
  &=\hsi((ab)^{-1})
  =\hsi(b^{-1}a^{-1})\\
  &=\hsi(a^{-1})+\hsi(b^{-1})+1\\
  &=-\hsi(a)-\hsi(b)+1\\
  &=-(\hsi(a)+\hsi(b)-1).
\end{align*}
To prove the forth property of $\hsi$ we consider the sequential set that corresponds to the standard basis $v$
and apply Lemma~\ref{lem-lift-glm}.
\end{proof}

We now formalize the properties of the function $\hsi$ in the following definition:

\begin{mydef}
\label{def-m-arf}
We denote by $\pi_1^0(P,p)$ the set of all non-trivial elements of $\pi_1(P,p)$
that can be represented by simple contours.
%that either do not belong to the kernel of the intersection form or are homologous to a hole or a puncture.
An {\it $m$-Arf function} is a function
$$\arf:\pi_1^0(P,p)\to\zm$$
satisfying the following conditions
\begin{enumerate}[1.]
\item
\label{arf-prop-conj}
$\arf(bab^{-1})=\arf(a)$ for any elements~$a,b\in\pi_1^0(P,p)$,
% such that the element~$bab^{-1}$ is in~$\pi_1^0(P,p)$,
\item
\label{arf-prop-inv}
$\arf(a^{-1})=-\arf(a)$ for any element~$a\in\pi_1^0(P,p)$,
\item
\label{arf-prop-cross}
$\arf(ab)=\arf(a)+\arf(b)$
for any
% elements~$a,b\in\pi_1^0(P,p)$ such that the element~$ab$ is in~$\pi_1^0(P,p)$ and the
elements~$a$ and~$b$ which can be represented by a pair of simple contours in $P$
intersecting in exactly one point~$p$ with $\<a,b\>\ne0$,
\item
\label{arf-prop-neg-seqset}
$\arf(ab)=\arf(a)+\arf(b)-1$
for any elements~$a,b\in\pi_1^0(P,p)$ such that the element~$ab$ is in~$\pi_1^0(P,p)$ and the elements~$a$ and~$b$ can be represented by a pair of simple contours in $P$
intersecting in exactly one point~$p$ with $\<a,b\>=0$
and placed in a neighbourhood of the point~$p$ as shown in Figure~\ref{fig-neg-pair}.
%\ie in such a way that the oriented contours $a$, $b$, and
%$(ab)^{-1}$ are freely homotopic to pairwise non-intersecting simple contours
%with orientation opposite to the one induced by the complex structure of the sphere with three
%holes that they cut out of $P$.
\end{enumerate}
% As for the function~$\hsi_{\lats}$, all equations involving~$\arf$ are to be understood as equations in~$\zm$.
\end{mydef}

\begin{rem}
One can prove that in the case $m=2$ there is a 1-1-correspondence between the $2$-Arf functions in the sense of Definition~\ref{def-m-arf}
and Arf functions in the sense of~\cite{Nbook}, Chapter~1, Section~7.
Namely, a function $\arf:\pi_1^0(P,p)\to\zd$ is a $2$-Arf function
if and only if $\om=1-\arf$ is an Arf function in the sense of~\cite{Nbook}.
\end{rem}

The following property of $m$-Arf functions follows immediately from Properties~\ref{arf-prop-neg-seqset} and~\ref{arf-prop-inv}
in Definition~\ref{def-m-arf}:

\begin{prop}
\label{arf-prop-seqset}
The equation $\arf(ab)=\arf(a)+\arf(b)+1$ is satisfied for any elements~$a,b\in\pi_1^0(P,p)$ such that
the element~$ab$ is in~$\pi_1^0(P,p)$ and the elements~$a$ and~$b$ can be represented by a pair of simple contours in $P$ intersecting in exactly one point~$p$ with $\<a,b\>=0$
and placed in a neighbourhood of the point~$p$ as shown in Figure~\ref{fig-pos-pair}.
\end{prop}

\begin{prop}
\label{arf-prop-liftexist}
For any standard basis
$$v=\{a_i,b_i~(i=1,\dots,g),c_i~(i=g+1,\dots,n)\}$$
of $\pi_1(P,p)$ we have
$$\sum\limits_{i=g+1}^n\,\arf(c_i)=(2-2g)-(n-g).$$
\end{prop}

\begin{proof}
We discuss the case $g=0$ first, and then we reduce the general case to the case $g=0$.
\begin{enumerate}[$\bullet$]
\item
Let $g=0$.
We prove that the statement is true for lifts of sequential sets of type $(0,l_h,l_p)$ by induction on~$l_h+l_p$.
In the case $l_h+l_p=3$ Proposition~\ref{arf-prop-seqset} implies
$$\arf(c_1c_2)=\arf(c_1)+\arf(c_2)+1,$$
while Property~\ref{arf-prop-inv} implies
$$\arf(c_1c_2)=\arf(c_3^{-1})=-\arf(c_3).$$
Combining these two equations we obtain
$$\arf(c_1)+\arf(c_2)+\arf(c_3)=-1.$$
Assume that the statement is true for $l_h+l_p\le n-1$ and consider the case $l_h+l_p=n$.
By our assumption
$$\arf(c_1\cdot c_2)+\arf(c_3)+\cdots+\arf(c_n)=2-(n-1)=(2-n)+1.$$
Moreover, by Proposition~\ref{arf-prop-seqset} we have
$$\arf(c_1c_2)=\arf(c_1)+\arf(c_2)+1,$$
The last two equations imply that
$$\sum\limits_{i=1}^n\,\arf(c_i)=2-n.$$
\item
We now consider the general case.
The set
$$(a_1,b_1a_1^{-1}b_1^{-1},\dots,a_g,b_ga_g^{-1}b_g^{-1},c_{g+1},\dots,c_n)$$
is a standard basis of a surface of type $(0,2g+l_h,l_p)$, hence
$$
  \sum\limits_{i=1}^g\,(\arf(a_i)+\arf(b_ia_i^{-1}b_i^{-1}))+\sum\limits_{i=g+1}^n\arf(c_i)
  =2-(n+g)
  =(2-2g)-(n-g).
$$
From Properties~\ref{arf-prop-conj} and~\ref{arf-prop-inv} of $m$-Arf functions we obtain that
$$
  \arf(b_i a_i^{-1}b_i^{-1})
  =\arf(a_i^{-1})
  =-\arf(a_i),
$$
and hence
$$\arf(a_i)+\arf(b_i a_i^{-1}b_i^{-1})=0.\qedhere$$
\end{enumerate}
\end{proof}

\begin{mydef}
\label{def-si-to-lats}
Let $\hsi_{\lats}:\pi_1(P,p)\to\zm$ be the function associated to a lift $\lats$
as in definition~\ref{def-hsi-to-lats},
then the function $\si_{\lats}:=\hsi_{\lats}|_{\pi_1^0(P,p)}$
is an $m$-Arf function according to
Lemma~\ref{lem-hsi-rules}, \ref{lem-inv}, and~\ref{lem-conj}.
We call the function $\si_{\lats}$ the {\it $m$-Arf function associated to the lift $\lats$.}
\end{mydef}

\subsection{Higher Arf functions and Dehn twists}

\label{arf-propert}

We recall the definition of the Dehn twists (\cite{Dehn}, \cite{Nbook} Chapter~1, Lemma~7.4).

\begin{mydef}
Under Dehn twists we understand some transformations from a standard basis
$$v=\{a_i,b_i~(i=1,\dots,g),c_i~(i=g+1,\dots,n)\}$$
of $\pi_1(P,p)$ to another standard basis
$$v'=\{a_i',b_i'~(i=1,\dots,g),c_i'~(i=g+1,\dots,n)\}.$$
Any of these transformations induces a homotopy class
of autohomeomorphisms of the surface~$P$,
which maps holes to holes and punctures to punctures.
By Theorem of Dehn~\cite{Dehn} the Dehn twists generate the whole group of such homotopy classes of autohomeomorphisms of~$P$.
There are Dehn twists of the following types:
\begin{align*}
1.\quad &a_1'=a_1b_1.\\
2.\quad &a_1'=(a_1a_2)a_1(a_1a_2)^{-1},\\
        &b_1'=(a_1a_2)a_1^{-1}b_1a_2^{-1}(a_1a_2)^{-1},\\
        &a_2'=a_1a_2a_1^{-1},\\
        &b_2'=b_2a_2^{-1}a_1^{-1}.\\
3.\quad &a_g'=(b_g^{-1}c_1)b_g^{-1}(b_g^{-1}c_1)^{-1},\\
        &b_g'=(b_g^{-1}c_1b_g)c_1^{-1}b_ga_gb_g^{-1}(b_g^{-1}c_1b_g)^{-1},\\
        &c_{g+1}'=b_g^{-1}c_{g+1}b_g.\\
4.\quad &a_k'=a_{k+1},\\
        &b_k'=b_{k+1},\\
        &a_{k+1}'=(c_{k+1}^{-1}c_k)a_k(c_{k+1}^{-1}c_k)^{-1},\\
        &b_{k+1}'=(c_{k+1}^{-1}c_k)b_k(c_{k+1}^{-1}c_k)^{-1}.\\
5.\quad &c_k'=c_{k+1},\\
        &c_{k+1}'=c_{k+1}^{-1}c_kc_{k+1}.
\end{align*}
Here $c_i=[a_i,b_i]$ for $i=1,\dots,g$,
in~4 we consider $k\in\{1,\dots,g\}$,
in~5 we consider $k\in\{g+1,\dots,n\}$.
If $a_i'$, $b_i'$ resp.\ $c_i'$ is not described explicitly, this means $a_i'=a_i$, $b_i'=b_i$ resp.\ $c_i'=c_i$.
\end{mydef}

Now we compute the values of $\arf$ on the standard basis~$v'$
from the values of $\arf$ on the standard basis~$v$.

\begin{lem}
\label{lem-dehn}
Let $\arf:\pi_1^0(P,p)\to\zm$ be an $m$-Arf function.
Let $D$ be a Dehn twist of the type described above.
Suppose that $D$ maps the standard basis
$$v=\{a_i,b_i~(i=1,\dots,g),c_i~(i=g+1,\dots,n)\}$$
into the standard basis
$$v'=D(v)=\{a_i',b_i'~(i=1,\dots,g),c_i'~(i=g+1,\dots,n)\}.$$
Let
$$\{\al_i,\be_i~(i=1,\dots,g),\ga_i~(i=g+1,\dots,n)\}$$
resp.
$$\{\al'_i,\be'_i~(i=1,\dots,g),\ga'_i~(i=g+1,\dots,n)\}$$
be the values of $\arf$ on the elements of $v$ resp.\ $v'$.
Then for the Dehn twists of types~1--5 we obtain
\begin{align*}
1.\quad &\al_1'=\al_1+\be_1.\\
2.\quad &\be_1'=\be_1-\al_1-\al_2-1,\\
        &\be_2'=\be_2-\al_2-\al_1-1.\\
3.\quad &\al_g'=-\be_g,\\
        &\be_g'=\al_g-\ga_{g+1}-1.\\
4.\quad &\al_k'=\al_{k+1},\\
        &\be_k'=\be_{k+1},\\
        &\al_{k+1}'=\al_k,\\
        &\be_{k+1}'=\be_k.\\
5.\quad &\ga_k'=\ga_{k+1},\\
        &\ga_{k+1}'=\ga_k.
\end{align*}
\end{lem}

\begin{proof}
We assume that the Dehn twist $D$ belongs to one of the types described in the definition above.
In the first case according to Property~\ref{arf-prop-cross} of $m$-Arf functions
%Corollary~\ref{cor-product-of-crossing-hyps}
we obtain
$$\arf(a_1')=\arf(a_1b_1)=\arf(a_1)+\arf(b_1).$$
In the second case according to Property~\ref{arf-prop-conj}
%Lemma~\ref{lem-conj}
we obtain
\begin{align*}
  \arf(a_1')&=\arf((a_1a_2)a_1(a_1a_2)^{-1})=\arf(a_1),\\
  \arf(b_1')&=\arf((a_1a_2)a_1^{-1}b_1a_2^{-1}(a_1a_2)^{-1})=\arf(a_1^{-1}b_1a_2^{-1}),\\
  \arf(a_2')&=\arf(a_1a_2a_1^{-1})=\arf(a_2),\\
  \arf(b_2')&=\arf(b_2a_2^{-1}a_1^{-1})
\end{align*}
According to Properties~\ref{arf-prop-cross} and~\ref{arf-prop-inv}
%Corollary~\ref{cor-product-of-crossing-hyps}
we obtain
\begin{align*}
  \arf(a_1^{-1}b_1)&=\arf(a_1^{-1})+\arf(b_1)=-\arf(a_1)+\arf(b_1),\\
  \arf(b_2a_2^{-1})&=\arf(b_2)+\arf(a_2^{-1})=\arf(b_2)-\arf(a_1).
\end{align*}

\myskip\noindent
In the following computations we illustrate the position of the contours on the surface
with figures showing the position of the axes of the corresponding elements in $\lat$.
Let
$$V=\{A_i, B_i~(i=1,\dots,g), C_i~(i=g+1,\dots,n)\}$$
be the sequential set corresponding to the standard basis~$v$.

\myskip\noindent
The mutual position of the axes of the elements $A_1^{-1}B_1$ and $A_2^{-1}$ is as in Figure~\ref{fig-axes-1}.

% Picture

\begin{figure}
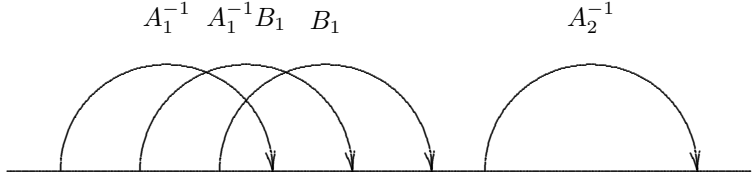

  \begin{center}
    \forcehmode
      \bgroup
        \beginpicture
          \setcoordinatesystem units <20 bp,20 bp>
          \multiput {\phantom{$\bullet$}} at -4 -1 10 3 /
          \circulararc 180 degrees from 1 0 center at -1 0
          \arrow <7pt> [0.2,0.5] from 0.995 0.04 to 1 0
          \put {$A_1^{-1}$} [b] <0pt,\baselineskip> at -1 2
          \circulararc 180 degrees from 4 0 center at 2 0
          \arrow <7pt> [0.2,0.5] from 3.995 0.04 to 4 0
          \put {$B_1$} [b] <0pt,\baselineskip> at 2 2
          \circulararc 180 degrees from 2.5 0 center at 0.5 0
          \arrow <7pt> [0.2,0.5] from 2.495 0.04 to 2.5 0
          \put {$A_1^{-1}B_1$} [b] <0pt,\baselineskip> at 0.5 2
          \circulararc 180 degrees from 9 0 center at 7 0
          \arrow <7pt> [0.2,0.5] from 8.995 0.04 to 9 0
          \put {$A_2^{-1}$} [b] <0pt,\baselineskip> at 7 2
          \plot -4 0 10 0 /
        \endpicture
      \egroup
  \end{center}
  \caption{Axes of $A_1^{-1}B_1$ and $A_2^{-1}$}
  \label{fig-axes-1}
\end{figure}

\noindent
According to Properties~\ref{arf-prop-neg-seqset} and~\ref{arf-prop-inv} we obtain
\begin{align*}
  \arf(b_1')
  &=\arf((a_1^{-1}b_1)a_2^{-1})\\
  &=\arf(a_1^{-1}b_1)+\arf(a_2^{-1})-1\\
  &=(\arf(b_1)-\arf(a_1))-\arf(a_2)-1\\
  &=\arf(b_1)-\arf(a_1)-\arf(a_2)-1.
\end{align*}
Similarly
$$\arf(b_2')=\arf(b_2)-\arf(a_2)-\arf(a_1)-1.$$
In the third case we obtain according to Properties~\ref{arf-prop-inv} and~\ref{arf-prop-conj}
% Lemma~\ref{lem-conj-of-hyps-and-pars}
\begin{align*}
  \arf(a_g')&=\arf((b_g^{-1}c_1)b_g^{-1}(b_g^{-1}c_1)^{-1})=\arf(b_g^{-1})=-\arf(b_g),\\
  \arf(b_g')&=\arf((b_g^{-1}c_1b_g)c_1^{-1}b_ga_gb_g^{-1}(b_g^{-1}c_1b_g)^{-1})=\arf(c_1^{-1}b_ga_gb_g^{-1}),\\
  \arf(c_1')&=\arf(b_g^{-1}c_1b_g)=\arf(c_1).
\end{align*}

% Picture

\begin{figure}
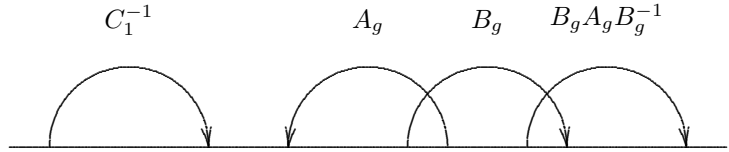

  \begin{center}
    \forcehmode
      \bgroup
        \beginpicture
          \setcoordinatesystem units <15 bp,15 bp>
          \multiput {\phantom{$\bullet$}} at -10 -1 8 3 /
          \circulararc 180 degrees from 1 0 center at -1 0
          \arrow <7pt> [0.2,0.5] from -2.995 0.04 to -3 0
          \put {$A_g$} [b] <0pt,\baselineskip> at -1 2
          \circulararc 180 degrees from 4 0 center at 2 0
          \arrow <7pt> [0.2,0.5] from 3.995 0.04 to 4 0
          \put {$B_g$} [b] <0pt,\baselineskip> at 2 2
          \circulararc 180 degrees from 7 0 center at 5 0
          \arrow <7pt> [0.2,0.5] from 6.995 0.04 to 7 0
          \put {$B_gA_gB_g^{-1}$} [b] <0pt,\baselineskip> at 5 2
          \circulararc 180 degrees from -5 0 center at -7 0
          \arrow <7pt> [0.2,0.5] from -5.005 0.04 to -5 0
          \put {$C_1^{-1}$} [b] <0pt,\baselineskip> at -7 2
          \plot -10 0 8 0 /
        \endpicture
      \egroup
  \end{center}
  \caption{Axes of $C_1^{-1}$ and $B_gA_gB_g^{-1}$}
  \label{fig-axes-2}
\end{figure}

\noindent
The mutual position of the axes of the elements $C_1^{-1}$ and $B_gA_gB_g^{-1}$ is as in Figure~\ref{fig-axes-2}.
According to Properties~\ref{arf-prop-neg-seqset} and~\ref{arf-prop-conj} we obtain
\begin{align*}
  \arf(b_g')
  &=\arf(c_1^{-1}\cdot(b_ga_gb_g^{-1}))\\
  &=\arf(c_1^{-1})+\arf(b_ga_gb_g^{-1})-1\\
  &=\arf(c_1^{-1})+\arf(a_g)-1.
\end{align*}

\noindent
In the forth and fifth case computations are easy, we only use Property~\ref{arf-prop-conj} of $m$-Arf functions.
\end{proof}

\subsection{Correspondence between higher Arf functions and higher spin structures}

\begin{lem}
\label{lem-arf-difference}
The difference $\arf_1-\arf_2:\pi_1^0(P,p)\to\zm$ of two Arf functions~$\arf_1$ and~$\arf_2$
induces a linear function $\ell:H_1(P;\zm)\to\zm$.
\end{lem}

\begin{proof}
Let
\begin{align*}
  v&=\{a_i,b_i~(i=1,\dots,g),c_i~(i=g+1,\dots,n)\}\\
  &=\{d_j~(j=1,\dots,n+g)\}
\end{align*}
be a standard basis of $\pi_1(P,p)$.
For an element $a\in\pi_1(P,p)$ let us denote by $[a]$ the image of~$a$ in $H_1(P;\z_m)$.
Let
\begin{align*}
  [v]&=\{[a_i],[b_i]~(i=1,\dots,g),[c_i]~(i=g+1,\dots,n)\}\\
  &=\{[d_j]~(j=1,\dots,n+g)\}
\end{align*}
be the induced basis of $H_1(P;\z_m)$.

\myskip
Let us define $\de:\pi_1^0(P,p)\to\zm$ by $\de=\arf_1-\arf_2$.
Properties of the Arf functions~$\arf_1$ and~$\arf_2$ imply
the following properties of the function~$\de$:
\begin{enumerate}[1.]
\item
\label{de-prop-conj}
$\de(bab^{-1})=\de(a)$,
\item
\label{de-prop-inv}
$\de(a^{-1})=-\de(a)$,
\item
\label{de-prop-cross}
$\de(ab)=\de(a)+\de(b)$
if the elements~$a$ and~$b$ can be represented by a pair of simple contours in $P$
intersecting in exactly one point~$p$
with $\<a,b\>\ne0$,
\item
\label{de-prop-neg-seqset}
$\de(ab)=\de(a)+\de(b)$
if the elements~$a$ and~$b$ can be represented by a pair of simple contours in $P$
intersecting in exactly one point~$p$
with $\<a,b\>=0$
and placed in a neighbourhood of the point~$p$ as shown in Figure~\ref{fig-neg-pair}.
\end{enumerate}

\myskip
Let us define a function $\ell:H_1(P;\z_m)\to\z_m$ on~$[v]$ by
$$\ell[d_j]=\de(d_j)$$
for all $j=1,\dots,n+g$
and then extend this function linearly on the whole $H_1(P;\z_m)$.
(Here we write $\ell[c]$ instead of $\ell([c])$ to simplify the notation.)
The linearity of $\ell$ implies the following properties:
\begin{enumerate}[1.]
\item
\label{ell-prop-conj}
$\ell[bab^{-1}]=\ell([b]+[a]-[b])=\ell[a]$,
\item
\label{ell-prop-inv}
$\ell[a^{-1}]=\ell(-[a])=-\ell[a]$,
\item
\label{ell-prop-cross}
$\ell[ab]=\ell[a]+\ell[b]$
if the elements~$a$ and~$b$ can be represented by a pair of simple contours in $P$
intersecting in exactly one point~$p$
with $\<a,b\>\ne0$,
\item
\label{ell-prop-neg-seqset}
$\ell[ab]=\ell[a]+\ell[b]$
if the elements~$a$ and~$b$ can be represented by a pair of simple contours in $P$
intersecting in exactly one point~$p$
with $\<a,b\>=0$
and placed in a neighbourhood of the point~$p$ as shown in Figure~\ref{fig-neg-pair}.
\end{enumerate}

\myskip\noindent
{\bf Claim:}
We claim that $\ell[a]=\de(a)$ for any element $a\in\pi_1^0(P,p)$.

\myskip\noindent
{\bf Observation:}
We first observe that if a standard basis
\begin{align*}
  v'&=\{a_i',b_i'~(i=1,\dots,g),c_i'~(i=g+1,\dots,n)\}\\
  &=\{d_j'~(j=1,\dots,n+g)\}
\end{align*}
is the image of the standard basis~$v$ under a Dehn twist,
then
$$\ell[d_j']=\de(d_j')$$
for all $j=1,\dots,n+g$.
Indeed, for the Dehn twists of type 1, 2, and 3 we obtain by Lemma~\ref{lem-dehn} and by linearity of~$\ell$
\begin{align*}
1.\quad \de(a_1')&=\de(a_1)+\de(b_1),\\
        \ell[a_1']&=\ell[a_1b_1]=\ell[a_1]+\ell[b_1],\\
2.\quad \de(b_1')&=\de(b_1)-\de(a_1)-\de(a_2),\\
        \ell[b_1']&=\ell[(a_1a_2)a_1^{-1}b_1a_2^{-1}(a_1a_2)^{-1}]=\ell[b_1]-\ell[a_1]-\ell[a_2],\\
        \de(b_2')&=\de(b_2)-\de(a_2)-\de(a_1),\\
        \ell[b_2']&=\ell[b_2a_2^{-1}a_1^{-1}]=\ell[b_2]-\ell[a_1]-\ell[a_2],\\
3.\quad \de(a_g')&=-\de(b_g),\\
        \ell[a_g']&=\ell[(b_g^{-1}c_1)b_g^{-1}(b_g^{-1}c_1)^{-1}]=-\ell[b_g],\\
        \de(b_g')&=\de(a_g)-\de(c_1),\\
        \ell[b_g']&=\ell[(b_g^{-1}c_1b_g)c_1^{-1}b_ga_gb_g^{-1}(b_g^{-1}c_1b_g)^{-1}]=\ell[a_g]-\ell[c_1].
\end{align*}
Here if $\de(d_j')$ and $\ell[d_j']$ are not described explicitly,
this means $\de(d_j')=\de(d_j)$ and $\ell[d_j']=\ell[d_j]$.

\myskip
Let $a$ be any contour not belonging to the kernel of the intersection form.
Then $a$ is an element $a_1'$ in some standard basis~$v'$ of $\pi_1(P,p)$.
By Theorem of Dehn~\cite{Dehn} the standard basis~$v'$ of $\pi_1(P,p)$
can be obtained from $v$ as an image under a sequence of Dehn twists.
Then the observation implies that $\de(a)=\ell[a]$.
Similarly any simple contour~$c$ homologous to a hole or a puncture
is an element $c'_{g+1}$ in some standard basis~$v'$ of $\pi_1(P,p)$,
hence $\de(c)=\ell[c]$.
Let us now consider any simple contour~$c$.
We can express the element~$c$ as a product of elements of some standard basis
\begin{align*}
  v'&=\{a_i',b_i'~(i=1,\dots,g),c_i'~(i=g+1,\dots,n)\}\\
  &=\{d_j'~(j=1,\dots,n+g)\}
\end{align*}
of $\pi_1(P,p)$.
Using the properties of the functions~$\de$ and~$\ell$
we compute $\de(c)$ as a linear combination of~$\de(d'_j)$
and $\ell[c]$ as the linear combination of~$\ell[d'_j]$ with the same coefficients.
For example for the contour~$c=c'_{g+1}c'_{g+2}$ we compute
\begin{align*}
  \de(c)&=\de(c'_{g+1})+\de(c'_{g+2}),\\
  \ell[c]&=\ell[c'_{g+1}]+\ell[c'_{g+2}].
\end{align*}
On the other hand according to the observation $\de(d'_j)=\ell[d'_j]$, hence $\de(c)=\ell[c]$.
So we have proved our claim.
Hence the function $\de=\arf_1-\arf_2$ induces a linear function on $H_1(P;\zm)$.

\myskip
On the contrary, if $\arf_0$ is an $m$-Arf function on~$P$ and $\ell$ is a linear function on $H_1(P;\zm)$,
% such that $$\sum\limits_{i=g+1}^n\,\ell[c_i]=0$$
then the function $\arf:\pi_1^0(P,p)\to\zm$ defined by
$$\arf(a)=\arf_0(a)+\ell[a]$$
is an $m$-Arf function on $P$.
Using the properties of the $m$-Arf function~$\arf_0$ and linearity of the function~$\ell$
we can prove that the function $\arf$ also has the properties of $m$-Arf functions.
For example, we show Property~\ref{arf-prop-neg-seqset}:
If $\<c_1,c_2\>=0$ and the elements~$c_1$ and~$c_2$ can be represented by a pair of simple contours in $P$
intersecting in exactly one point~$p$ and placed in a neighbourhood of the point~$p$ as is shown in Figure~\ref{fig-neg-pair},
then
$$\arf(c_1c_2)=\arf(c_1)+\arf(c_2)-1.$$
Linearity of the function~$\ell$ and Property~\ref{arf-prop-neg-seqset} for $\arf_0$ imply
\begin{align*}
  \arf(c_1c_2)
  &=\arf_0(c_1c_2)+\ell[c_1c_2]\\
  &=(\arf_0(c_1)+\arf_0(c_2)-1)+(\ell[c_1]+\ell[c_2])\\
  &=\arf(c_1)+\arf(c_2)-1.
\end{align*}
The proofs of the other properties are similar.
\end{proof}

\begin{cor}
\label{cor-arf-affine}
The set $\Arf^{P,m}$ of all $m$-Arf functions on $\pi_1^0(P,p)$ has a structure of an affine space,
\ie the set $\{\arf-\arf_0\st\arf\in\Arf^{P,m}\}$ is a free module over~$\zm$ for any $\arf_0\in\Arf^{P,m}$.
\end{cor}

\begin{cor}
\label{cor-arf-def-by-gen}
An $m$-Arf function is uniquely determined by its values on the elements of some standard basis of $\pi_1(P,p)$.
\end{cor}

\begin{cor}
\label{cor-number-of-arfs}
Let $P$ be a surface of type $(g,l_h,l_p)$.
Then the number of different $m$-Arf functions on $\pi_1^0(P,p)$ is equal to
\begin{enumerate}[$\bullet$]
\item
$m^{2g+l_h+l_p-1}$, if $l_h+l_p>0$,
\item
$m^{2g+l_h+l_p}$, if $l_h+l_p=0$ and $2g-2=0~\mod~m$,
\item
0, otherwise.
\end{enumerate}
\end{cor}

\begin{thm}
\label{thm-corresp}
There is a 1-1-correspondence between
\begin{enumerate}[1)]
\item $m$-spin structures on~$P=\hyp/\lat$.
% \ie complex line bundles over~$P$ such that $L^{\otimes m}\cong TP$.
\item lifts of~$\lat$ into $\Gm$,
\item $m$-Arf functions $\arf:\pi_1^0(P,p)\to\zm$.
\end{enumerate}
\end{thm}

\begin{proof}
According to Proposition~\ref{spinstructures-lifts}
there is a 1-1-correspondence between $m$-spin structures on~$P$
and the lifts of $\lat$ into $\Gm$.
In Definition~\ref{def-si-to-lats} we attached to any lift $\lats$ of $\lat$
into $\Gm$ an $m$-Arf function $\si_{\lats}$.
On the other hand we can attach to any $m$-Arf function $\arf$ a subset of~$\Gm$
$$\lats_{\arf}:=\{g\in\Gm\st\pi(g)\in\lat,~\levm(g)=\arf(\Phi(\pi(g)))\},$$
where $\pi:\Gm\to G$ is the covering map.
It remains to prove that this subset of~$\Gm$ is actually a lift of~$\lat$.
Let
\begin{align*}
  v&=\{a_i, b_i~(i=1,\dots,g), c_i~(i=g+1,\dots,n)\}\\
  &=\{d_j~(j=1,\dots,n+g)\}
\end{align*}
be a standard basis of $\pi_1(P,p)$
and let
\begin{align*}
  \bV
  &=\{\bA_i, \bB_i~(i=1,\dots,g), \bC_i~(i=g+1,\dots,n)\}\\
  &=\{\bD_j~(j=1,\dots,n+g)\}\\
  &=\{\Phi^{-1}(d_j)~(j=1,\dots,n+g)\}
\end{align*}
be the corresponding sequential set.
Let
\begin{align*}
  V
  &=\{A_i, B_i~(i=1,\dots,g), C_i~(i=g+1,\dots,n)\}\\
  &=\{D_j~(j=1,\dots,n+g)\}
\end{align*}
be a lift of the sequential set~$\bV$, \ie $\pi(D_j)=\bD_j$, such that $\levm(D_j)=\arf(d_j)$.
Then we obtain
$$
  \sum\limits_{i=g+1}^n\,\levm(C_i)
  =\sum\limits_{i=g+1}^n\,\arf(c_i)
  =(2-2g)-(n-g),
$$
hence by Lemma~\ref{lem-lift-glm} we obtain
$$\prod\limits_{i=1}^g\,[A_i,B_i]\cdot\prod\limits_{i=g+1}^n\,C_i=e.$$
This implies according to Lemma~\ref{lem-liftonerel} that
the subgroup $\lats$ of $\Gm$ generated by $V$ is a lift of $\lat$ into $\Gm$.
Let us compare the corresponding Arf function $\arf_{\lats}$ with the Arf function $\arf$.
We have
$$\arf_{\lats}(d_j)=\levm(D_j)=\arf(d_j)$$
for all~$j$ \ie the Arf functions $\arf_{\lats}$ and $\arf$ coincide on the standard basis~$v$.
Thus by Corollary~\ref{cor-arf-def-by-gen} the Arf functions $\arf_{\lats}$ and $\arf$ coincide on the whole $\pi_1^0(P,p)$.
From the definition of $\arf_{\lats}$ and $\lats_{\arf}$ we see that this implies that $\lats=\lats_{\arf}$,
hence $\lats_{\arf}$ is indeed a lift of~$\lat$ into~$\Gm$.
It is clear from the definitions that the mappings $\lats\mapsto\arf_{\lats}$ and $\arf\mapsto\lats_{\arf}$
are inverse to each other.
\end{proof}

\begin{rem}
Let us define an {\it $\infty$-Arf function} as a function $\arf:\pi_1^0(P,p)\to\z$ satisfying
the four properties in Definition of $m$-Arf functions (Definition~\ref{def-m-arf}).
Then slightly modifying the previous discussion of $m$-Arf functions we obtain
similar results for $\infty$-Arf functions:
\begin{enumerate}[1)]
\item
The set $\Arf^{P,\infty}$ of all $\infty$-Arf functions on $\pi_1^0(P,p)$ has a structure of an affine space
associated to $H^1(P;\z)$.
\item
For a surface $P=\hyp/\lat$ there is a 1-1-correspondence
between $\infty$-Arf functions $\arf:\pi_1^0(P,p)\to\z$
and lifts of the group~$\lat$ into the universal cover $\tG$ of $\PSL$.
\end{enumerate}
\end{rem}

\begin{cor}
The set of all $m$-spin structures on~$P$ has a structure of an affine space.
\end{cor}

\begin{cor}
Let $P$ be a surface of type $(g,l_h,l_p)$.
Then the number of different $m$-spin structures on $P$ is equal to
\begin{enumerate}[$\bullet$]
\item
$m^{2g+l_h+l_p-1}$, if $l_h+l_p>0$,
\item
$m^{2g+l_h+l_p}$, if $l_h+l_p=0$ and $2g-2=0~\mod~m$,
\item
0, otherwise.
\end{enumerate}
\end{cor}

\section{Moduli spaces of higher spin structures}

\label{sec-moduli-spaces}

We study the moduli space of $m$-spin hyperbolic Riemann surfaces, i.e.\ hyperbolic Riemann surfaces with an $m$-spin structure.
Using Proposition~\ref{spinstructures-lifts}, we define the moduli space of $m$-spin structures as the space of conjugacy classes of subgroups~$\lats$ in~$\Gm$
such that the restriction of the covering map~$\Gm\to G=\PSL$ to $\lats$ is an isomorphism between $\lats$ and a Fuchsian group~$\lat$.

\bigskip\noindent
The projection~$\lats\mapsto\lat$ from the moduli space of Riemann surfaces with an $m$-spin structure to the moduli space of Riemann surfaces is a finite ramified covering.

\myskip
\subsection{Topological classification of higher Arf functions}

\myskip
There is a 1-1-cor\-res\-pon\-dence (see Theorem~\ref{thm-corresp})
between $m$-spin structures on a surface~$P$ and $m$-Arf functions on~$\pi_1^0(P,p)$.
This correspondence allows us to reduce the problem of finding the number of connected components
of the moduli space of higher spin Riemann surfaces
to the problem of finding the number of orbits of the action of the group of autohomeomorphisms
on the set of $m$-Arf functions.
We describe the orbit of an $m$-Arf function under the action of the group of homotopy classes of surface autohomeomorphisms.

\myskip
Let $P$ be a Riemann surface of type $(g,l_h,l_p)$ with $l_h+l_p=n-g$.
Let $p\in P$.

\begin{mydef}
\label{def-arf-inv}
Let $\arf:\pi_1^0(P,p)\to\zm$ be an $m$-Arf function.
For~$g=1$ we define the {\it Arf invariant} $\de=\de(P,\arf)$ as
$$\de=\gcd(m,\arf(a_1),\arf(b_1),\arf(c_2)+1,\dots,\arf(c_n)+1),$$
where
$$\{a_1,b_1,c_i~(i=2,\dots,n)\}$$
is a standard basis of the fundamental group $\pi_1(P,p)$.
(The Arf invariant~$\de$ is well defined as it is invariant under the transformations described in Lemma~\ref{lem-dehn}.)
For~$g>1$ and even $m$ we define the {\it Arf invariant} $\de=\de(P,\arf)$ as $\de=0$
if there is a standard basis
$$\{a_i,b_i~(i=1,\dots,g),c_i~(i=g+1,\dots,n)\}$$
of the fundamental group $\pi_1(P,p)$ such that
$$\sum\limits_{i=1}^g(1-\arf(a_i))(1-\arf(b_i))=0~\mod~2$$
and as $\de=1$ otherwise.
For~$g>1$ and odd $m$ we set $\de=0$.
\end{mydef}

%It remains to show that the the Arf invariant is well defined, \ie that the
%value of the sum in the definition does not dependent on the choice of the
%standard basis. To this end we show using the explicit description of the
%action of Dehn twists on $m$-Arf functions in Lemma~\ref{}
%we can prove that the parity of the sum
%$$\sum\limits_{i=1}^g(1-\arf(a_i))(1-\arf(b_i))$$
%is preserved under the Dehn twists,
%hence the Arf invariant $\de$ is well defined.

\begin{mydef}
\label{def-type-hyp-arf-func}
Let $\arf:\pi_1^0(P,p)\to\zm$ be an $m$-Arf function
and
$$v=\{a_i,b_i~(i=1,\dots,g),c_i~(i=g+1,\dots,n)\}$$
a standard basis of the fundamental group $\pi_1(P,p)$.
In the standard basis~$v$ the elements $c_{g+1},\dots,c_{g+l_h}$ correspond to the holes
and the elements $c_{g+l_h+1},\dots,c_n$ correspond to the punctures.
The set of holes
$$c_{g+1},\dots,c_{g+l_h}$$
is divided into $m$ sets
$$D_j=\{c_i\st\arf(c_i)=j,~g<i\le g+l_h\}.$$
We denote by $\num^h_j=\num^h_j(P,\arf)$, $j\in\zm$, the cardinality of the set $D_j$.
Similarly, the punctures
$$c_{g+l_h+1},\dots,c_n$$
are divided into $m$ sets of the form
$$E_j=\{c_i\st\arf(c_i)=j,~g+l_h<i\le n\}.$$
Let $\num^p_j=\num^p_j(P,\arf)$, $j\in\zm$, be the cardinality of the set $E_j$.
We have
$$\sum_{j\in\zm}\,\num^h_j=l_h\quad\text{and}\quad\sum_{j\in\zm}\,\num^p_j=l_p.$$
By the {\it type of the $m$-Arf function}~$(P,\arf)$ we mean the tuple
$$(g,\de,\num^h_0,\dots,\num^h_{m-1},\num^p_0,\dots,\num^p_{m-1}),$$
where $\de$ is the Arf invariant of $\arf$ defined above.
\end{mydef}

\begin{lem}
\label{lem-normalize}
Assume that~$g>1$.
Let $\arf:\pi_1^0(P,p)\to\zm$ be an $m$-Arf function,
then there is a standard basis
$$v=\{a_i, b_i~(i=1,\dots,g), c_i~(i=g+1,\dots,n)\}$$
of $\pi_1(P,p)$ such that
$$(\arf(a_1),\arf(b_1),\dots,\arf(a_g),\arf(b_g))=(0,\xi,1,\dots,1)$$
with $\xi\in\{0,1\}$ and
$$
  \arf(c_{g+1})\le\dots\le\arf(c_{g+l_h}),
  \quad
  \arf(c_{g+l_h+1})\le\dots\le\arf(c_n).
$$
If $m$ is odd or there is a contour around a hole or a puncture such that
the value of $\arf$ on this contour is even, then the basis can be chosen in
such a way that $\xi=1$.
\end{lem}

\begin{proof}

\myskip\noindent
\begin{enumerate}[a)]
\item
\label{proof-part-reformulation}
Let us fix some standard basis
$$v_0=\{a_i, b_i~(i=1,\dots,g), c_i~(i=g+1,\dots,n)\}$$
and consider the sequence of values of the Arf function $\arf$ on this basis~$v_0$
\begin{align*}
  &(\al_1, \be_1,\dots,\al_g,\be_g,\ga_{g+1},\dots,\ga_n)\\
  &=(\arf(a_1),\arf(b_1),\dots,\arf(a_g),\arf(b_g),\arf(c_{g+1}),\dots,\arf(c_n)).
\end{align*}
If $m$ is even and there is a hole or a puncture, such that for the corresponding contour~$c$ the value~$\arf(c)$ is even,
then we can choose the standard basis~$v_0$ in such a way, that the value~$\arf(c_{g+1})$ is even.
Any other standard basis $v$ is an image of the basis $v_0$ under an autohomeomorphisms of the surface,
\ie under a sequence of Dehn twist.
Hence according to Lemma~\ref{lem-dehn} the sequence of values of $\arf$ on the basis~$v$
is the image of the corresponding sequence with respect to the basis $v_0$ under
the group generated by the transformations that change
the first $2g$ components $(\al_1,\be_1,\dots,\al_g,\be_g)$ as follows
\begin{align*}
  &1a.~(\al_1,\be_1,\dots,\al_i,\be_i,\dots,\al_g,\be_g)\\
  &\qquad\mapsto(\al_1,\be_1,\dots,\al_i\pm\be_i,\be_i,\dots,\al_g,\be_g),\\
  &1b.~(\al_1,\be_1,\dots,\al_i,\be_i,\dots,\al_g,\be_g)\\
  &\qquad\mapsto(\al_1,\be_1,\dots,\al_i,\be_i\pm\al_i,\dots,\al_g,\be_g),\\
  &2.~(\al_1,\be_1,\dots,\al_i,\be_i,\dots,\al_j,\be_j,\dots,\al_g,\be_g)\\
  &\qquad\mapsto(\al_1,\be_1,\dots,\al_i,\be_i-\al_j-1,\dots,\al_j,\be_j-\al_i-1,\dots,\al_g,\be_g),\\
  &3.~(\al_1,\be_1,\dots,\al_{g-1},\be_{g-1},\al_g,\be_g)\\
  &\qquad\mapsto(\al_1,\be_1,\dots,\al_{g-1},\be_{g-1},-\be_g,\al_g-\ga_{g+1}-1),\\
  &4.~(\al_1,\be_1,\dots,\al_i,\be_i,\dots,\al_j,\be_j,\dots,\al_g,\be_g)\\
  &\qquad\mapsto(\al_1,\be_1,\dots,\al_j,\be_j,\dots,\al_i,\be_i,\dots,\al_g,\be_g),\\
  &5a.~(\al_1,\be_1,\dots,\al_i,\be_i,\dots,\al_g,\be_g)\\
  &\qquad\mapsto(\al_1,\be_1,\dots,-\al_i,-\be_i,\dots,\al_g,\be_g),\\
  &5b.~(\al_1,\be_1,\dots,\al_i,\be_i,\dots,\al_g,\be_g)\\
  &\qquad\mapsto(\al_1,\be_1,\dots,-\be_i,\al_i,\dots,\al_g,\be_g)
\end{align*}
and change $(\ga_{g+1},\dots,\ga_n)$ by all possible permutations of $(\ga_{g+1},\dots,\ga_{g+l_h})$ and $(\ga_{g+l_h+1},\dots,\ga_n)$.
The inequalities between the values of $\arf$ on the elements $c_i$ are easy
to satisfy, because the transformation group contains all possible permutations of $(\ga_{g+1},\dots,\ga_{g+l_h})$
and $(\ga_{g+l_h+1},\dots,\ga_n)$.
Our aim is to show that, if $m$ is odd or one of the numbers $\ga_{g+1},\dots,\ga_n$ is even,
any tuple $(\al_1,\be_1,\dots,\al_g,\be_g)$ can be transformed into the $2g$-tuple
$$(0,1,1,\dots,1),$$
while otherwise any such tuple can be transformed into one of the tuples
$$(0,0,1,\dots,1)\quad\hbox{or}\quad(0,1,1,\dots,1).$$
\item
\label{proof-part-minus2}
We claim that the group of transformations described in~(\ref{proof-part-reformulation})
contains the transformation of the form
$$
  (\dots,\al_i,\be_i,\dots,\al_j,\be_j,\dots)
  \mapsto
  (\dots,\al_i,\be_i-2,\dots,\al_j,\be_j,\dots):
$$
Let us assume $(i,j)=(1,2)$ in order to simplify the notation.
We apply transformations~2, 5a, again 2 and again 5a
and obtain
\begin{align*}
  &(\al_1,\be_1,\al_2,\be_2,\dots)\\
  &\mapsto(\al_1,\be_1-\al_2-1,\al_2,\be_2-\al_1-1,\dots)\\
  &\mapsto(\al_1,\be_1-\al_2-1,\al_2,\be_2-\al_1-1,\dots)\\
  &\mapsto(\al_1,\be_1-\al_2-1,-\al_2,-(\be_2-\al_1-1),\dots)\\
  &\mapsto(\al_1,\be_1-\al_2-1-(-\al_2)-1,-\al_2,-\be_2+\al_1+1-\al_1-1,\dots)\\
  &=(\al_1,\be_1-2,-\al_2,-\be_2,\dots)\\
  &\mapsto(\al_1,\be_1-2,\al_2,\be_2,\dots).
\end{align*}
\item
\label{proof-part-zanuli}
Furthermore, we claim that the group of transformations contains a transformation of the form
$$
  (\al_1,\be_1,\al_2,\be_2,\dots)\mapsto(0,\xi,1,1,\dots),
  \quad\hbox{where}\quad\xi\in\{0,1\}:
$$
With the help of the transformation described in~(\ref{proof-part-minus2}) we can transform
$$
  (\al_1,\be_1,\al_2,\be_2,\dots)\mapsto(\al_1',\be_1',\al_2',\be_2',\dots),
  \quad\hbox{where}\quad\al_1',\be_1',\al_2',\be_2'\in\{0,1\}.
$$
Applying the inverse of transformation 2 we obtain
$$(0,0,0,0,\dots)\mapsto(0,1,0,1,\dots).$$
Applying transformation 1a resp.~1b we obtain
$$(0,1,\dots)\mapsto(1,1,\dots)\quad\hbox{and}\quad(1,0,\dots)\mapsto(1,1,\dots)$$

\item
By successive application of the transformation described in~(\ref{proof-part-zanuli})
we obtain that any $2g$-tuple can be transformed to a tuple of the form $(0,\xi,1,\dots,1)$.

\item
If $m=2r+1$ is odd, then we use the transformation
described in~(\ref{proof-part-minus2}) to map
\begin{align*}
  (0,1,1,\dots,1)
  &\mapsto(0,1-2\cdot(r+1),1,\dots,1)\\
  &=(0,-m,1,\dots,1)=(0,0,1,\dots,1),
\end{align*}
hence the $2g$-tuple $(0,0,1,\dots,1)$ can be transformed
into the tuple $(0,1,1,\dots,1)$.

\item
If $m$ and $\ga_{g+1}=2r$ are even, then we use
the transformation~4, the transformation~3,
successive application of the transformation
described in~(\ref{proof-part-minus2}),
and the transformation~4 again to map
%$$(0,0,1,\dots,1)\mapsto(0,0-2r-1,1,\dots,1)\mapsto(0,1,1,\dots,1),$$
\begin{align*}
  (0,0,1,\dots,1)
  &\mapsto(1,\dots,1,0,0)\\
  &\mapsto(1,\dots,1,0,0-2r-1)\\
  &\mapsto(1,\dots,1,0,1),\\
  &\mapsto(0,1,1,\dots,1),
\end{align*}
hence also in this case the $2g$-tuple $(0,0,1,\dots,1)$
can be transformed into the tuple $(0,1,1,\dots,1)$.\qedhere
\end{enumerate}
\end{proof}

\begin{lem}
\label{lem-normalize-genus1}
Assume that~$g=1$.
Let $\arf:\pi_1^0(P,p)\to\zm$ be an $m$-Arf function, then there is a standard basis
$$v=\{a_1, b_1, c_i~(i=2,\dots,n)\}$$
of $\pi_1(P,p)$ such that
$$(\arf(a_1),\arf(b_1))=(\de,0),$$
where $\de=\gcd(m,\al_1,\be_1,\ga_2+1,\dots,\ga_n+1)$, and
$$
  \arf(c_2)\le\dots\le\arf(c_{l_h+1}),
  \quad
  \arf(c_{l_h+2})\le\dots\le\arf(c_n).
$$
\end{lem}

\begin{proof}
\label{proof-normalize-genus1}
Let us fix some standard basis
$$v_0=\{a_1, b_1, c_i~(i=2,\dots,n)\}$$
and consider the sequence of values of the Arf function $\arf$ on this basis~$v_0$
$$(\al_1, \be_1,\ga_2,\dots,\ga_n)=(\arf(a_1),\arf(b_1),\arf(c_2),\dots,\arf(c_n)).$$
Any other standard basis $v$ is an image of the basis $v_0$ under an autohomeomorphisms of the surface,
\ie under a sequence of Dehn twist.
Hence according to Lemma~\ref{lem-dehn} the sequence of values of $\arf$ on the basis~$v$
is the image of the corresponding sequence with respect to the basis $v_0$ under
the group generated by the transformations that change the first two components $(\al_1,\be_1)$ as follows
\begin{align*}
  &1a.~(\al_1,\be_1)\mapsto(\al_1\pm\be_1,\be_1),\\
  &1b.~(\al_1,\be_1)\mapsto(\al_1,\be_1\pm\al_1),\\
  &3.~(\al_1,\be_1)\mapsto(-\be_1,\al_1-\ga_2-1),\\
  &5a.~(\al_1,\be_1)\mapsto(-\al_1,-\be_1),\\
  &5b.~(\al_1,\be_1)\mapsto(-\be_1,\al_1).
\end{align*}
and change $(\ga_2,\dots,\ga_n)$ by all possible permutations of
$$(\ga_2,\dots,\ga_{1+l_h})\quad\text{and}\quad(\ga_{1+l_h+1},\dots,\ga_n).$$
The inequalities between the values of $\arf$ on the elements $c_i$ are easy
to satisfy, because the transformation group contains all possible permutations of $(\ga_2,\dots,\ga_{1+l_h})$ and $(\ga_{1+l_h+1},\dots,\ga_n)$.
Using transformations~1a, 1b, 5a, 5b any tuple~$(\al_1,\be_1)$ can be transformed into the tuple~$(\gcd(m,\al_1,\be_1),0)$.
Using the transformation~3 any tuple~$(\al_1,\be_1)$ can be transformed into the tuple~$(\gcd(m,\al_1,\be_1,\ga_2+1,\dots,\ga_n+1),0)$.
%On the other hand, $\gcd(m,\al_1,\be_1,\ga_2+1,\dots,\ga_n+1)$ can not be changed by transformations~1a, 1b, 3, 5a, 5b.
\end{proof}

\begin{rem}
%The group of autohomeomorphisms of a surface $P$ generates a group $\calA$
% of automorphisms of all $m$-spin structures on~$P$.
An autohomeomorphism of a surface $P$ induces an automorphism of $m$-spin structures on~$P$.
Let $\calA$ be the corresponding group of such automorphisms of $m$-spin structures on~$P$.
On the other hand, any autohomeomorphism generates an element of~$\MathOpSp(2g,\z)$, where $g$ is the genus of~$P$.
Lemmas~\ref{proof-normalize-genus1}, \ref{def-type-hyp-arf-func} imply that for two autohomeomorphisms
the corresponding elements in~$\calA$ differ if the corresponding elements in~$\MathOpSp(2g,\z_m)$ differ.
Thus we obtain a homomorphism~$f:\calA\to\MathOpSp(2g,\z_m)$.
Using Dehn twists of types~1,2 and~4 we can show that $f$ is an epimorphism.
Lemmas~\ref{proof-normalize-genus1}, \ref{def-type-hyp-arf-func} imply that $\ker(f)=d\cdot T$,
where $T$ is the group of all parallel translations on the affine space of all $m$-spin structures.
Using Dehn twists of types~1--5 and Lemmas~\ref{proof-normalize-genus1}, \ref{def-type-hyp-arf-func}
we are able to determine the number~$d$.
If~$g>1$ then~$d=2$ if $m$ is even and all $n$ are odd, and $d=1$ otherwise.
If~$g=1$ then $d=\gcd(m,\al_1,\be_1,\ga_2+1,\dots,\ga_n+1)$.
An algebraic-geometric proof of this affirmation is contained in~\cite{Jarvis:2000}~(Lemma 3.3.6).
\end{rem}

\begin{thm}
\label{thm-is-type}
A tuple $t=(g,\de,\num^h_0,\dots,\num^h_{m-1},\num^p_0,\dots,\num^p_{m-1})$, where
$$\sum_{j\in\zm}\,\num^h_j=l_h\quad\text{and}\quad\sum_{j\in\zm}\,\num^p_j=l_p,$$
is the type of a hyperbolic $m$-Arf function on a Riemann surface of type
$(g,l_h,l_p)$
if and only if it has the following properties:
\begin{enumerate}[(a)]
\item
If $g>1$ and $m$ is odd, then $\de=0$.
\item
If $g>1$ and $m$ is even and $\num^h_j+\num^p_j\ne0$ for some even~$j\in\zm$, then $\de=0$.
\item
If $g=1$ then $\de$ is a divisor of~$m$ and~$\gcd(\{j+1\st\num^h_j+\num^p_j\ne0\})$.
\item
The following degree condition is satisfied
$$\sum\limits_{j\in\zm}\,j\cdot(\num^h_j+\num^p_j)=(2-2g)-(l_h+l_p).$$
\end{enumerate}
\end{thm}

\begin{proof}
Let us first assume that the tuple $t$ is the type of a hyperbolic $m$-Arf function
on a Riemann surface of type $(g,l_h,l_p)$ and that $\arf$ is such a function.
Then Proposition~\ref{arf-prop-liftexist} implies that
$$
  (2-2g)-(n-g)
  =\sum\limits_{i=g+1}^n\,\arf(c_i)
  =\sum\limits_{j\in\zm}\,j\cdot(\num^h_j+\num^p_j).
$$

\myskip\noindent
Let us assume that~$g>1$. If $m$ is odd or ($m$ is even and $\num^h_j+\num^p_j\ne0$ for some even~$j\in\zm$),
then according to Lemma~\ref{lem-normalize} there is a standard basis
$$v=\{a_i, b_i~(i=1,\dots,g), c_i~(i=g+1,\dots,n)\}$$
of $\pi_1(P,p)$ such that
$$(\arf(a_1),\arf(b_1),\dots,\arf(a_g),\arf(b_g))=(0,1,1,\dots,1),$$
hence $\de(P,\arf)=0$ by definition.
Now let us assume that $t=(g,\de,\num^h_j,\num^p_j)$ is a tuple satisfying the conditions~(1) and (2).
Let us fix some standard basis
$$v_0=\{a_i, b_i~(i=1,\dots,g), c_i~(i=g+1,\dots,n)\}$$
The condition~(1) together with Proposition~\ref{arf-prop-liftexist} and Corollary~\ref{cor-arf-def-by-gen} imply
that there exist Arf functions $\arf^0$ and $\arf^1$
such that $\num^h_j(\arf^{\de})=\num^h_j$, $\num^p_j(\arf^{\de})=\num^p_j$ for $j\in\zm$ and
$$(\arf^{\de}(a_1),\arf^{\de}(b_1),\dots,\arf^{\de}(a_g),\arf^{\de}(b_g))=(0,1-\de,1,\dots,1).$$
The equation $\de(\arf^0)=0$ is satisfied by definition.
It remains to prove, in the case that $m$ is even and all $\ga_i$ are odd, that $\de(\arf^1)=1$.
To this end we observe, using the explicit description of the action of Dehn twists on $m$-Arf functions in Lemma~\ref{lem-dehn},
that the parity of the sum
$$\sum\limits_{i=1}^g(1-\arf(a_i))(1-\arf(b_i))$$
is preserved under the Dehn twists and hence is equal to~$1$ modulo~$2$ for any standard basis.

\myskip\noindent
In the case~$g=1$ the proof is similar. In this case we observe that
$$\gcd(m,\arf(a_1),\arf(b_1),\arf(c_2)+1,\dots,\arf(c_n)+1)$$
is preserved under the Dehn twists.
\end{proof}

\subsection{Connected components of the moduli space}

\label{components}

We recall the results on the moduli space of Riemann surfaces from~\cite{Nbook}.
\myskip
Let $\Ga_{g,n}$ be the group generated by the elements
$$v=\{a_1,b_1,\dots,a_g,b_g,c_{g+1},\dots,c_n\}$$
with a single defining relation
$$\prod\limits_{i=1}^g\,[a_i,b_i]\prod\limits_{i=g+1}^n\,c_i=1.$$

\myskip
We denote by $\tTgll$ the set of monomorphisms $\psi:\Ga_{g,n}\to\Aut(\hyp)$
such that
$$\psi(v)=\{a_i^{\psi}, b_i^{\psi}~(i=1,\dots,g), c_i^{\psi}~(i=g+1,\dots,n)\}$$
is a sequential set of type $(g,l_h,l_p)$.
Here $l_h+l_p=n-g$ and we assume that $6g+3l_h+2l_p>6$.

\myskip
The group $\Aut(\hyp)$ acts on $\tTgll$ by conjugation.
We set
$$\Tgll=\tTgll/\Aut(\hyp).$$

\myskip
We parametrise the space $\tTgll$ by the fixed points and shift parameters of the elements of the sequential sets $\psi(v)$.
We use here the following result similar to Theorem of Fricke and Klein~\cite{FK}.

\begin{thm} {\rm (\cite{Nbook} (Chapter~1, Theorem~4.1))}

\myskip\noindent
The space $\Tgll$ is diffeomorphic to an open domain in $\r^{6g+3l_h+2l_p-6}$,
which is homeomorphic to $\r^{6g+3l_h+2l_p-6}$.
\end{thm}

For an element $\psi:\Ga_{g,n}\to\Aut(\hyp)$ of $\tTgll$ we write
$$\tMod^{\psi}=\tMod_{g,l_h,l_p}^{\psi}=\{\al\in\Aut(\Ga_{g,n})\st\psi\circ\al\in\tTgll\}.$$
One can show that $\tMod^{\psi}$ does not depend on~$\psi$, hence we write $\tMod$ instead of~$\tMod^{\psi}$.
Let $I\tMod$ be the subgroup of all inner automorphisms of $\Ga_{g,n}$ and let
$$\Modgll=\Mod=\tMod/I\tMod.$$
We now recall the description of the moduli space of Riemann surfaces

\begin{thm} {\rm (\cite{Nbook}, Chapter~1, Section~5)}
\label{moduli-of-surfaces}

\myskip\noindent
The group $\Mod=\Modgll$ and the group of homotopy classes of orientation preserving autohomeomorphisms of the surface of type $(g,l_h,l_p)$
are naturally isomorphic.
The group $\Modgll$ acts naturally on $\Tgll$ by diffeomorphisms.
This action is discrete.
The quotient set $\Tgll/\Modgll$ can be identified naturally with the moduli space $\Mgll$ of Riemann surfaces of type $(g,l_h,l_p)$.
\end{thm}

\begin{mydef}
We denote by
$$\Smgll(t)=\Smgll(g,\de,\num^h_0,\dots,\num^h_{m-1},\num^p_0,\dots,\num^p_{m-1})$$
the set of all $m$-spin structures on all Riemann surfaces of type~$(g,l_h,l_p)$
such that the associated $m$-Arf function is of type
$$t=(g,\de,\num^h_0,\dots,\num^h_{m-1},\num^p_0,\dots,\num^p_{m-1}).$$
\end{mydef}

\begin{thm}
\label{top-type-comp}
Let $t=(g,\de,\num^h_0,\dots,\num^h_{m-1},\num^p_0,\dots,\num^p_{m-1})$ be a
tuple that satisfies the conditions of Theorem~\ref{thm-is-type},
\ie the space $\Smgll(t)$ is not empty.
Then the space $\Smgll(t)$ is homeomorphic to $\Tgll/\Modmgll(t)$, where $\Tgll$ is
homeomorphic to~$\r^{6g+3l_h+2l_p-6}$
and $\Modmgll(t)$ acts on $\Tgll$ as a subgroup of finite index in the group $\Modgll$.
\end{thm}

\begin{proof}
Let us consider an element~$\psi$ of the space $\Tgll$.
By definition $\psi$ is an homomorphism $\psi:\Ga_{g,n}\to\Aut(\hyp)$.
To the homomorphism $\psi$ we attach a Riemann surface $P_{\psi}=\hyp/\psi(\Ga_{g,n})$,
a standard basis
$$v_{\psi}=\{a^{\psi}_i, b^{\psi}_i~(i=1,\dots,g), c^{\psi}_i~(i=g+1,\dots,n)\}$$
of $\pi_1(P_{\psi},p)$,
and an $m$-Arf function $\arf_{\psi}$ on this surface given by
\begin{align*}
  &(\arf_{\psi}(a_1^{\psi}),\arf_{\psi}(b_1^{\psi}))=(\de,0)\quad\text{if}\quad g=1,\\
  &(\arf_{\psi}(a_1^{\psi}),\arf_{\psi}(b_1^{\psi}),\arf_{\psi}(a_2^{\psi}),\arf_{\psi}(b_2^{\psi}),\dots,\arf_{\psi}(a_g^{\psi}),\arf_{\psi}(b_g^{\psi}))\\
  &
  \begin{aligned}
  &=(0,1-\de,1,\dots,1)\quad\text{if}\quad g>1,\\
  (\arf_{\psi}(c_{g+1}^{\psi}),\dots,\arf_{\psi}(c_{g+l_h}^{\psi}))
  &=(
     \underbrace{0,\dots,0}_{\num^h_0},
     \underbrace{1,\dots,1}_{\num^h_1},
     \dots,
     \underbrace{m-1,\dots,m-1}_{\num^h_{m-1}}
    ),\\
  (\arf_{\psi}(c_{g+l_h}^{\psi}),\dots,\arf_{\psi}(c_n^{\psi}))
  &=(
     \underbrace{0,\dots,0}_{\num^p_0},
     \underbrace{1,\dots,1}_{\num^p_1},
     \dots,
     \underbrace{m-1,\dots,m-1}_{\num^p_{m-1}}
    ).
  \end{aligned}
\end{align*}
By Theorem~\ref{thm-corresp}, the $m$-Arf function $\arf_{\psi}$ on the surface $P_{\psi}$
corresponds to an $m$-spin structure $\de_{\psi}$ on~$P_{\psi}$.
The correspondence $\psi\mapsto \de_{\psi}$ defines a map
$$\Tgll\to\Smgll(t).$$
According to Theorem~\ref{thm-is-type} this map is surjective.
Let $\Modmgll(t)$ be the subgroup of $\Aut(P_{\psi})=\Modgll$ that preserves the $m$-Arf function $\arf_{\psi}$.
For any point in $\Smgll(t)$ its preimage in $\Tgll$ consists of an orbit of the subgroup~$\Modmgll(t)$.
Thus
$$\Smgll(t)=\Tgll/\Modmgll(t).\qedhere$$
%where $\Modmgll(t)$ is a discrete group.
\end{proof}

\bigskip\noindent
Summarizing the results of Theorems~\ref{thm-is-type} and \ref{top-type-comp} we obtain the following

\bigskip
\begin{thm*}

\noindent
\begin{enumerate}[1)]
\item
Two $m$-spin structures are in the same connected component
of the space of all $m$-spin structures on hyperbolic Riemann surfaces
iff they are of the same type.
In other words, the connected components of the space of all $m$-spin
structures are those sets $\Smgll(t)$ that are not empty.
\item
The set $\Smgll(t)$ is not empty iff
$t=(g,\de,\num^h_0,\dots,\num^h_{m-1},\num^p_0,\dots,\num^p_{m-1})$
has the following properties:
\begin{enumerate}[(a)]
\item
If $g>1$ and $m$ is odd, then $\de=0$.
\item
If $g>1$ and $m$ is even and $\num^h_j+\num^p_j\ne0$ for some even~$j\in\zm$, then $\de=0$.
\item
If $g=1$ then $\de$ is a divisor of~$m$ and~$\gcd(\{j+1\st\num^h_j+\num^p_j\ne0\})$.
\item
The following degree condition is satisfied
$$\sum\limits_{j\in\zm}\,j\cdot(\num^h_j+\num^p_j)=(2-2g)-(l_h+l_p).$$
\end{enumerate}
\item
Any connected component~$\Smgll(t)$ of the space of all $m$-spin structures on hyperbolic Riemann surfaces of type~$(g,l_h,l_p)$
is homeomorphic to
$$\r^{6g+3l_h+2l_p-6}/\Modmgll(t),$$
where $\Modmgll(t)$ is a subgroup of finite index in the group~$\Modgll$
and acts discrete on $\r^{6g+3l_h+2l_p-6}$.
\end{enumerate}
\end{thm*}

%\begin{rem}
%Number of possible types of $m$-Arf functions: XXX
%\end{rem}

%\nocite{*}

\bibliographystyle{amsalpha}
\bibliography{arf}

\end{document}